# Properties, Proved and Conjectured, of Keller, Mycielski, and Queen Graphs

Witold Jarnicki, Krakow, Poland; witoldj@gmail.com
Wendy Myrvold, University of Victoria, Victoria, British Columbia, Canada; wendym@cs.uvic.ca
Peter Saltzman, Berkeley, California, USA; saltzman.pwa@gmail.com
Stan Wagon, Macalester College, St. Paul, Minnesota, USA; wagon@macalester.edu

**Abstract.** We prove several results about three families of graphs. For queen graphs, defined from the usual moves of a chess queen, we find the edge-chromatic number in almost all cases. In the unproved case, we have a conjecture supported by a vast amount of computation, which involved the development of a new edge-coloring algorithm. The conjecture is that the edge-chromatic number is the maximum degree, except when simple arithmetic forces the edge-chromatic number to be one greater than the maximum degree. For Mycielski graphs, we strengthen an old result that the graphs are Hamiltonian by showing that they are Hamilton-connected (except $M_3$, which is a cycle). For Keller graphs $G_d$, we establish, in all cases, the exact value of the chromatic number, the edge-chromatic number, and the independence number; and we get the clique covering number in all cases except $5 \le d \le 7$. We also investigate Hamiltonian decompositions of Keller graphs, obtaining them up to $G_6$.

## 1. Introduction

Inspired by computational experiments, we prove several results about some families of graphs. We show in §5 that all Mycielski graphs (except $M_3$, which is a 5-cycle) are Hamilton-connected. In §6, we establish the size of a maximum independent set for all Keller graphs and investigate some other parameters, determining the chromatic number of both the graphs and their complements, and also the edge-chromatic number. In particular, we prove that the edge-chromatic number of each Keller graph equals its degree. We also find the clique covering number for all cases except dimension 5, 6, and 7. And in Sections 2–4 we present a detailed study of queen graphs, resolving the edge-chromatic number in most cases.

Recall that the problem of coloring the edges of a graph is much simpler than the classic vertex-coloring problem. There are only two possibilities for the edge-chromatic number because of Vizing's classic theorem [BM, sec. 6.2] that the edge-chromatic number $\chi'(G)$ is either $\Delta(G)$ or $\Delta(G)+1$, where $\Delta(G)$ is the maximum vertex degree; the first case is called *class 1*; the second, *class 2*. Let $n_e(G)$ denote the number of edges and $\rho(G)$ the number of edges in a maximum matching. Some graph have too many edges to be class 1. An *overfull* graph $G$ is one for which $n_E(G) > \Delta(G) \left\lfloor \frac{n_v(G)}{2} \right\rfloor$. For such a graph, $\Delta(G)\rho(G) < n_e(G)$, and this inequality implies that $G$ must be class 2; so any overfull graph is class 2. The reason for this is that each color class is a matching and so has size at most $\rho(G)$; if class 1, the number of colored edges would be at most $\Delta(G)\,\rho(G)$ which is too small to capture all edges. In Section 2, we present results and an intriguing conjecture related to edge coloring of the standard queen graph $Q_{m,n}$: the conjecture is that $Q_{m,n}$ is class 1 whenever it is not overfull. Computation and proofs yield the truth of this conjecture for $m \le 10$ and all values of $n \ge m$; the exact conjecture is that the queen is class 1 for $n \le \frac{1}{3}\left(2m^3 - 11m + 12\right)$. In Theorem 7, we prove this for $n \le \frac{1}{2}\left(m^2 - 3m + 2\right)$. For the extensive computations, we developed a general edge-coloring algorithm that succeeded in finding class-1 colorings for some queen graphs having over two million edges.

Our notation is fairly standard: $K_n$ is the complete graph on $n$ vertices; $C_n$ is an $n$-cycle; $n_v(G)$ is the number of vertices (or the *order*) of $G$; $\chi(G)$ is the chromatic number; $\chi_{\text{frac}}(G)$ is the fractional chromatic number; $\alpha(G)$ is the size of a largest independent set; $\omega(G)$ is the size of a largest clique; $\theta(G)$ is the clique covering number (same as $\chi(G^c)$). Occasionally $G$ will be omitted from these functions where the context is clear. A vertex of $G$ is called *major* if its degree equals $\Delta(G)$. Graphs are always simple graphs, with the exception of some queen graph discussions, where multigraphs appear.

A *Hamiltonian path* (resp. *cycle*) is a path (resp. cycle) that passes through all vertices and does not intersect itself. A graph is *Hamiltonian* if there is a Hamiltonian cycle; a graph is *Hamilton-connected* (HC) if, for any pair $u$, $v$ of vertices, there is a Hamilton path from the $u$ to $v$. We will make use of Fournier's Theorem [F1, F2] states that a graph is class 1 if the subgraph induced by the vertices of maximum degree is a forest.

We thank Joan Hutchinson for a careful reading and helpful suggestions, and David Pike for the interesting comment about edge coloring pioneer F. Walecki.

## 2. Rook and Bishop Graphs

The family of (not necessarily square) queen graphs presents a number of well-known combinatorial challenges. In this and the following two sections, we study the Vizing classification of queen graphs, a problem that turns out to have unexpected complexity. Queens on a chessboard can make all the moves of rooks and bishops, and thus queen graphs are the union of the rook graph and (white and black) bishop graphs. We therefore start, in Section 2, by looking at rook and bishop graphs separately. Then in Sections 3 and 4, we will show how these rook and bishop results lead to a variety of class-1 queen colorings.

It is well known that rook graphs behave similarly to their one-dimensional cousins, the complete graphs: they are class 2 if and only if both dimensions are odd. Perhaps more surprising, all bishop graphs are class 1. These two results already suffice to show that queen graphs are class 1 when at least one of the dimensions is even: just take the union of a class 1 rook coloring and a class 1 bishop coloring. When both dimensions are odd, however, the classification of queen graphs becomes much harder. A straightforward counting argument shows that such odd queen graphs are eventually class 2: for $m$ and $n$ odd, $Q_{m,n}$ is class 2 if $n \geq \frac{1}{3}\left(2m^3 - 11m + 18\right)$. On the other hand, we prove below (Thm. 7) that for $m$ and $n$ odd, $Q_{m,n}$ is class 1 if $m \leq n \leq \frac{1}{2}\left(m^2 - 3m + 2\right)$. As we will also show, the method we use cannot produce class-1 colorings all the way up to the cubic limit, and thus we leave essentially open the problem of determining whether there are any class-2 queen graphs when $\frac{1}{2}\left(m^2 - 3m + 4\right) \leq n \leq \frac{1}{3}\left(2m^3 - 11m + 12\right)$. We do, however, describe an algorithmic approach that gives lots of data to support the conjecture that there are no such graphs.

Recall that bishops move diagonally on a chessboard and rooks (Fig. 1) move horizontally or vertically. Because a queen can move diagonally, horizontally, or vertically, $Q_{m,n}$, the graph of queen moves, is the union of its two edge subgraphs $B_{m,n}$ and $R_{m,n}$, where $B_{m,n}$ denotes the graph of bishop moves on an $m \times n$ board, and $R_{m,n}$ denotes the rook graph; the latter is just the Cartesian product $K_m \,\square\, K_n$. The bishop graph is disconnected: it is the union of graphs corresponding to a white bishop and a black bishop (where we take the lower left square as being white). We will use $WB_{m,n}$ for the white bishop graph. It is natural to try to get edge-coloring results for the queen by combining such results for bishops and rooks, so we review the situation for those two pieces. The classic result on edge-coloring complete graphs is also essential, so we start there. Lucas [L, p. 177] attributes the first part of Proposition 1 to Felix Walecki.

**Proposition 1.** $\chi'(K_n)$ is $n-1$ when $n$ is even (and so the graph is class 1) and $n$ when $n$ is odd (the graph is class 2). Moreover, when $n$ is odd every coloring has the property that no missing color at a vertex is repeated. Also, for all $n$, if $M$ is a maximum matching of $K_n$, then $\chi'(K_n \backslash M) = n-1$.

**Proof.** For the even case, take the vertices to be $v_i$ where $v_1, \ldots, v_{n-1}$ are the vertices of a regular $(n-1)$-gon, and $v_n$ is the center. Use color $i$ on $v_i \leftrightarrow v_n$ and on edges perpendicular to this edge. For $n$ odd, one can use a regular $n$-gon to locate all the vertices and use $n$ colors for the exterior $n$-cycle; then color any other edge with the same color used for the exterior edge that parallels it. (Alternatively, add a dummy vertex $v_{n+1}$ and use the even-order result, discarding at the end any edges involving the dummy vertex.) Note that $K_n$, with $n$ odd, is overfull, so the preceding coloring is optimal. Further, the coloring has the property that the missing colors at the vertices are $1, 2, \ldots, n$. This phenomenon, that no missing color is repeated, is easily seen to hold for any class-2 coloring of $K_n$, with $n$ odd. Because $\chi'(K_n) = n_e(K_n)/\rho(K_n)$, each color class in any optimal coloring of $K_n$ is a maximum matching. These graphs are edge transitive, so all maximum matchings are the same, which yields the final assertion of the Proposition. □

**Theorem 2.** The rook graph $R_{m,n}$ is class 1 except when both dimensions are odd, in which case it is overfull, and so is class 2.

**Proof.** For the class-2 result, we have $\Delta = m + n - 2$, and $n_e = m\binom{n}{2} + n\binom{m}{2}$, which leads to $n_e - \frac{1}{2}(m\,n - 1)\,\Delta = \frac{\Delta}{2} > 0$. For class 1, the even case is trivial by Proposition 1, since we can use colors 1 through $n - 1$ on each row and $n$ through $m + n - 2$ on each column. For the case of $m$ even and $n$ odd (which suffices by symmetry), use colors 1 through $n$ on each complete row, and ensure that color 1 is missing at the vertices in the first column, color 2 is missing on the second column, and so on. The color set consisting of $i$, $n + 1, \ldots, n + m - 2$ can be used on the $i$th column. □

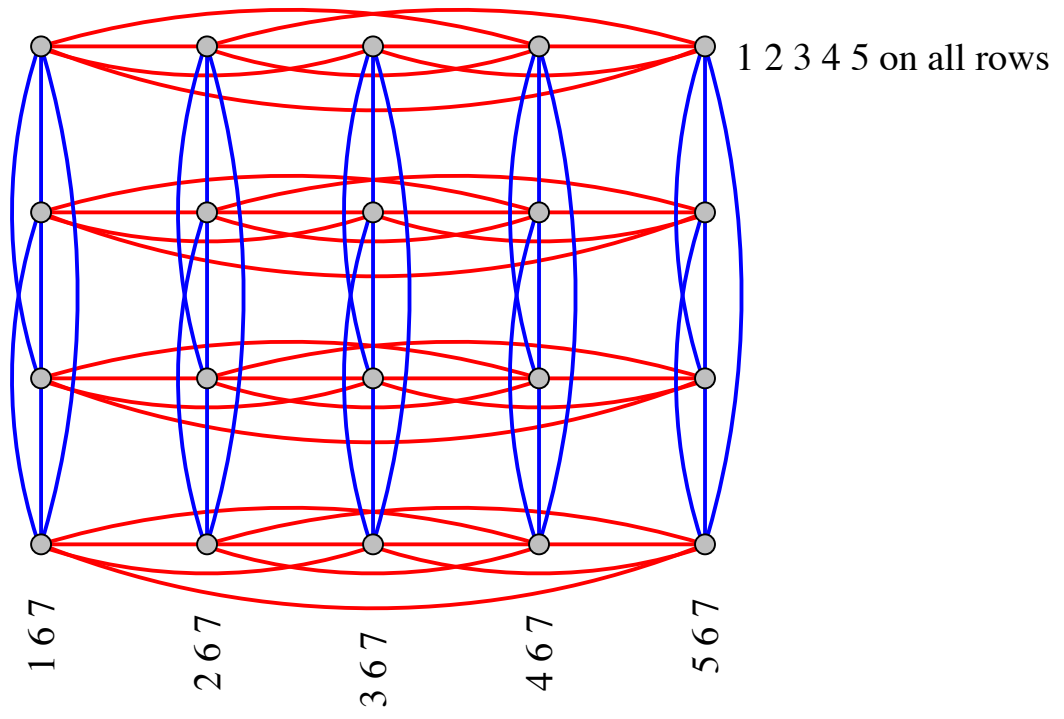

Figure 1. The rook graph $R_{4,5}$ is class 1. Use colors 1–5 on the rows, with color $i$ missing on the vertices in the $i$th column. Then use $i$, 6 and 7 on the $i$th column.

For the class-2 case of the preceding result one can easily give an explicit $\Delta + 1$ coloring, either by the method used in Theorem 6 or the ladder method of Proposition 8.

The bishop situation had not been investigated until the work of Saltzman and Wagon, who proved that all bishop graphs are class 1.

**Theorem 3 [SW, SW1].** All bishop graphs $B_{m,n}$ are class 1.

**Proof.** Assume $m \le n$. We have $\Delta(B_{m,n}) = 2m - 2$, except for one case: if $m$ is even, then $\Delta(B_{m,m}) = 2m - 3$. The graph can be decomposed into paths as follows. Note that any bishop edge is a diagonal line with a natural "length": the Euclidean distance between the vertices divided by $\sqrt{2}$. Let $G_1^+$ consist of all edges of length 1 having negative slope and paths of length $m - 1$, with edges having positive slope (in Fig. 2 this graph is the set of green and red edges). This subgraph consists of disjoint paths; the edges of each path can be 2-colored. Define $G_1^-$ the same way, but with the slopes reversed. Get the full family by defining $G_1^+, G_1^-, G_2^+, G_2^-, \ldots, G_{\lfloor m/2 \rfloor}^+, G_{\lfloor m/2 \rfloor}^-$, where $G_i^\pm$ is defined similarly to $G_1^\pm$, but using edges of length $i$ and $m - i$. The proof that these edge subgraphs partition the bishop edges is easy (see [SW]). Each of these subgraphs, being a collection of disjoint paths, can be 2-edge colored (for definiteness and because it plays a role in later work, we will always use the first of the two colors on the leftmost edge of each path; and the resulting coloring will be referred to as the *canonical bishop edge-coloring*). When $m$ is odd the color count is $4\,\frac{m-1}{2} = 2m - 2$. When $m$ is even and $n > m$, the graphs $G_{m/2}^+$ and $G_{m/2}^-$ coincide and the color count is $4\left(\frac{m}{2} - 1\right) + 2 = 2m - 2$. But when $m$ is even and $n = m$, then $G_{m/2}^+ = G_{m/2}^-$ and this subgraph consists of only disjoint edges; it is therefore 1-colorable and the color count is $4\left(\frac{m}{2} - 1\right) + 1 = 2m - 3$. Note that for odd $m$, some of the edge subgraphs when restricted to the black bishop will be empty, but that is irrelevant. The black bishop will use fewer colors, but the colors are disjoint from the ones used for the white bishop and it is the latter that determines $\chi'$. □

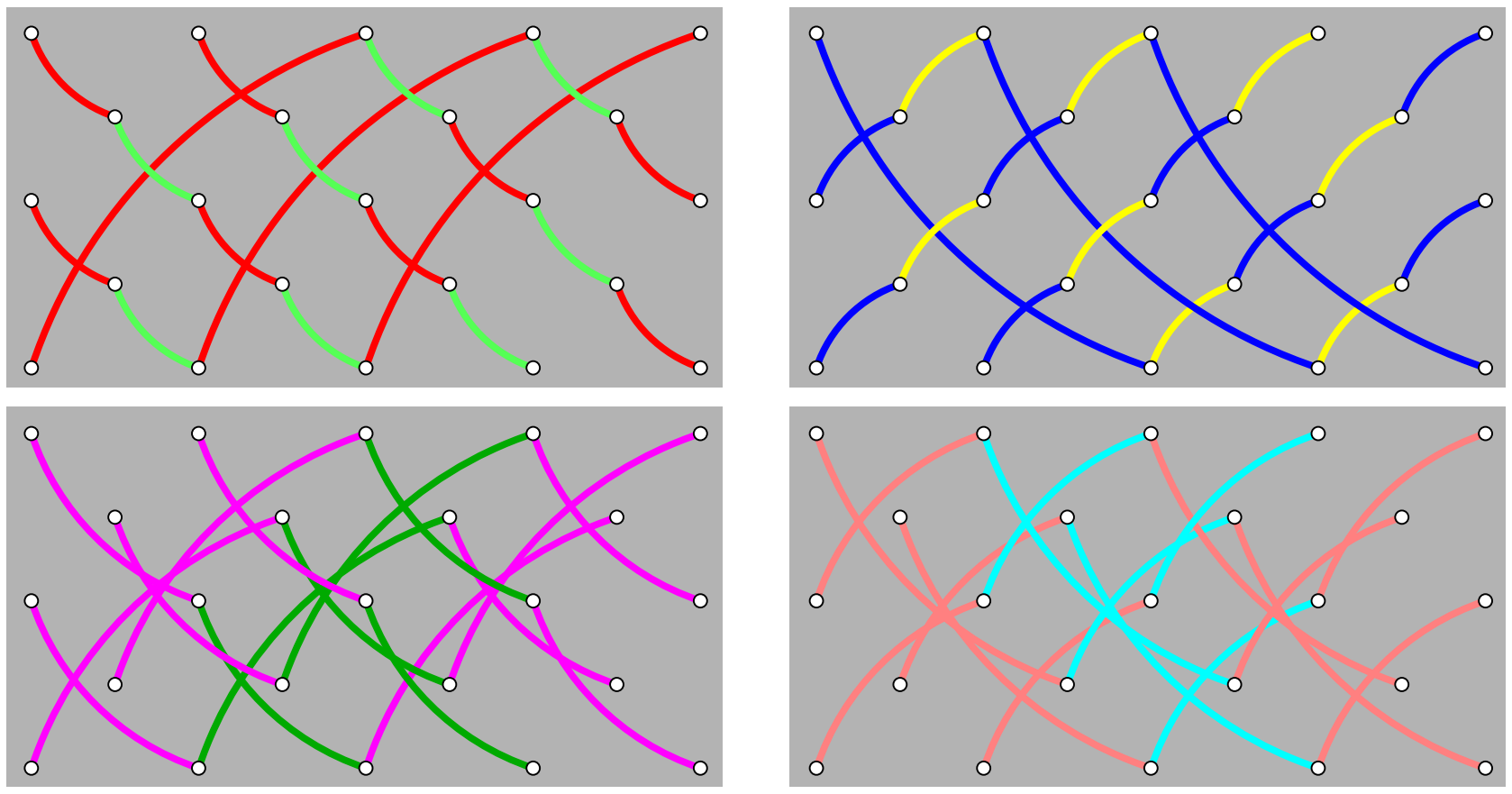

Figure 2. Top: The red and green edges form the subgraph $G_1^+$ of $WB_{5,9}$, the white bishop graph; blue and yellow form $G_1^-$. Bottom: The purple and green edges form $G_2^+$; cyan and pink are $G_2^-$. The total color count is 8, the maximum degree of the graph.

An interesting and useful type of coloring is one in which one color is as rare as can be. The next lemma shows that, for $B_{n,n}$ with $n$ odd, the canonical coloring is such that the rarest color occurs the smallest possible number of times: once.

**Lemma 4.** In the canonical coloring of $B_{n,n}$, $n$ odd, the rarest color occurs on one edge only.

**Proof.** Referring to the subgraphs of Theorem 3's proof, all the paths in the black bishop part of $G^-_{(m-1)/2}$ are isolated edges, and the same is true for the white bishop except for the single path $Z \longleftrightarrow X \longleftrightarrow Y$ where $X$ is the central vertex and $Y$, $Z$ are, respectively, the upper-right and lower-left corners (Fig. 3). Therefore the coloring used in the proof of Theorem 3 will use the last color only on $X \longleftrightarrow Y$. $\square$

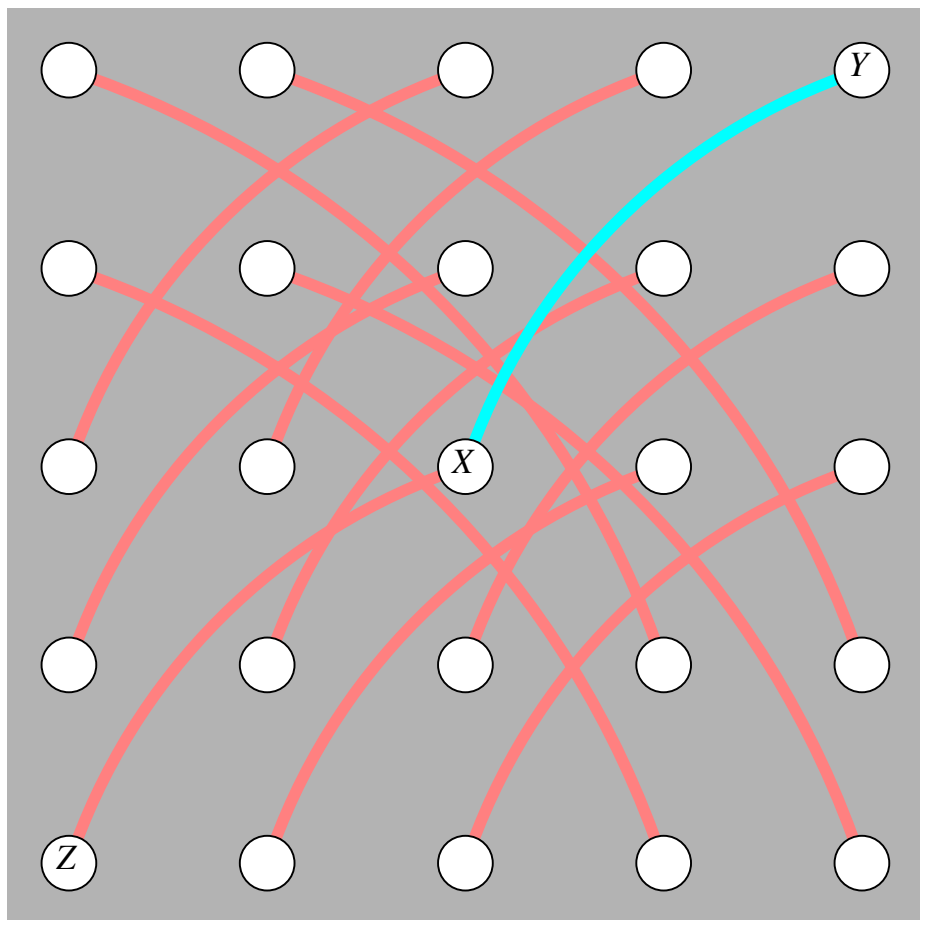

Figure 3. The edges shown form the subgraph $G_3^-$ of $B_{5,5}$. This subgraph consists of many isolated edges and one 2-edge path: $Z \longleftrightarrow X \longleftrightarrow Y$, and so can be edge-colored so that one color occurs only once, at $X \longleftrightarrow Y$.

A version of Lemma 4, with proof similar to the one given, but requiring some color switching, holds for even bishops; we do not need the result so just sketch the proof.

**Lemma 5.** The bishop graph $B_{2k,2k}$ admits a class-1 coloring where one of the colors appears on only 2 edges; and the 2 cannot be replaced by 1.

**Sketch of proof.** The proof focuses on the white bishop and uses some color switching in the subgraphs $G_1^-$ and $G_{k-1}^+$ to get the uniquely appearing color at the edge connecting the two major vertices. That 2 is best possible follows from the fact that the graph splits into two isomorphic components. □

In fact, this whole discussion can be generalized. Let $M(G)$ be the number of major vertices of $G$. Then in any coloring of $G$ using $\Delta(G)$ colors, each color must occur at least $\lceil \frac{1}{2} M \rceil$ times (because all colors appear at every major vertex). Moreover, a coloring using the rarest color exactly $\lceil \frac{1}{2} M \rceil$ times cannot exist unless there is a perfect matching (or almost perfect matching if $M(G)$ is odd) of the major subgraph, since such a matching is needed to make each account for two major vertices. Call a class-1 coloring of $G$ *extremal* if the rarest color occurs exactly $\lceil \frac{1}{2} M \rceil$ times. Then Lemma 4 states that $B_{n,n}$ admits an extremal class-1 coloring when $n$ is odd, and Lemma 5 implies that $B_{n,n}$ does not admit such a coloring when $n$ is even. Computations support the following conjecture, where $WB$ denotes the white bishop graph.

**Conjecture 1.** $WB_{m,n}$ always admits an extremal class-1 coloring.

## 3. Queen Graphs

The graph of queen moves on an $m \times n$ chessboard is the *queen graph* $Q_{m,n}$ (Fig. 4 shows $Q_{3,3}$). The vertices of $Q_{m,n}$ are arranged in an $m \times n$ grid and each vertex is adjacent to all vertices in the same row, in the same column, and on the same diagonal or back diagonal. We always assume $m \leq n$. Easy counting and summation leads to

$$\Delta(Q_{m,n}) = \begin{cases} 3m + n - 5 & \text{if } m = n \text{ and } n \text{ is even} \\ 3m + n - 4 & \text{otherwise} \end{cases}$$

$$n_e(Q_{m,n}) = \tfrac{1}{6} m\left(2 - 2m^2 - 12n + 9mn + 3n^2\right).$$

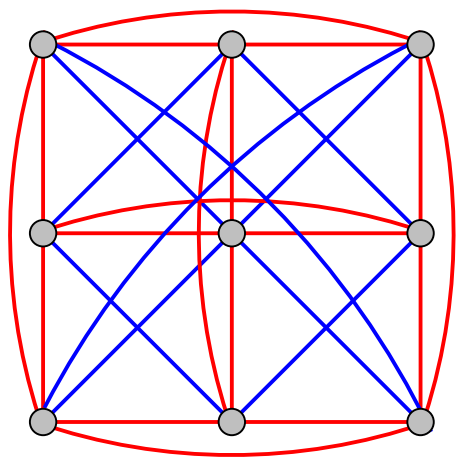

Figure 4. The queen graph $Q_{3,3}$ has 28 edges and maximum degree 8. It combines the bishop graph (blue) with the rook graph (red).

Queen graphs have served as a challenging benchmark for vertex coloring algorithms [OEIS]. The values of $\chi(Q_{n,n})$ are known for $n \leq 25$. For $11 \leq n \leq 25$, $\chi(Q_{n,n}) = n$; the first open case is $Q_{26,26}$ for which $\chi$ is known to be either 26, 27, or 28. So a case involving $26^2 = 676$ vertices is unresolved. As discussed later in this section, we developed an algorithm that succeeds in finding the edge-chromatic number for cases as large as $Q_{11,707}$, which has 7777 vertices and requires coloring almost three million edges. But many cases are resolved by relatively straightforward arguments and Theorems 6 and 7 find the edge-chromatic number in all cases except when $m$, $n$ are odd and $\frac{1}{2}\left(m^2 - 3m + 4\right) \leq n \leq \frac{1}{3}\left(2\, m^3 - 11\, m + 12\right)$.

Any queen graph is the edge-union of a bishop subgraph and a rook subgraph (Fig. 4); the maximum degrees add: $\Delta(Q_{m,n}) = \Delta(B_{m,n}) + \Delta(R_{m,n})$. Theorem 3 shows that all bishop graphs are class 1 and Theorem 2 shows that the rooks are class 1 except when both $m$ and $n$ are odd. Thus one can often get a class-1 queen coloring by forming the union of optimal colorings of the bishop and rook subgraphs. When the rook is class 1, one can simply combine class-1 colorings for the rook and bishop to get a class-1 queen coloring. This yields the class-1 part of the next theorem in all cases except one: $Q_{n,n}$, $n$ odd.

**Theorem 6.** (Joseph DeVincentis, Witold Jarnicki, and Stan Wagon) The queen graph $Q_{m,n}$ is class 1 if at least one of $m$ and $n$ is even, or if $m$ and $n$ are equal and odd. The graph is class 2 if $m$ and $n$ are odd and $n \geq \frac{1}{3}\left(2m^3 - 11m + 18\right)$.

**Proof.** The last assertion follows from the fact that $Q_{m,n}$ in that case is overfull and therefore is class 2. This is because the overfull condition becomes

$$\frac{1}{2}(mn-1)\,(3m+n-4) \le \frac{1}{6}\,m\left(2-2m^2-12\,n+9mn+3n^2\right)-1,$$

which simplifies to the stated inequality $n \ge \frac{1}{3}\left(2m^3-11m+18\right)$.

The class-1 result is proved by combining class-1 colorings of the rook and bishop subgraphs, except in the one case that the rooks are class 2. Thus a different argument is needed for $Q_{n,n}$ where $n$ is odd.

Consider $Q_{n,n}$ with $n$ odd. The central vertex is the ony vertex of maximum degree, so the result follows from Fournier's theorem (Section 1). It also follows from Theorem 7 below, but we can give a direct construction of a class-1 coloring, using a special property of the square bishop graph.

Start with a coloring of $B_{n,n}$ as in Lemma 4. Then we can color the corresponding rook graph $R_{n,n}$ using only new colors in such a way that a color is free to replace color $2n-2$ at its single use on the bishop edge $X \leftrightarrow Z$. The result will be a class-1 coloring of $Q_{n,n}$.

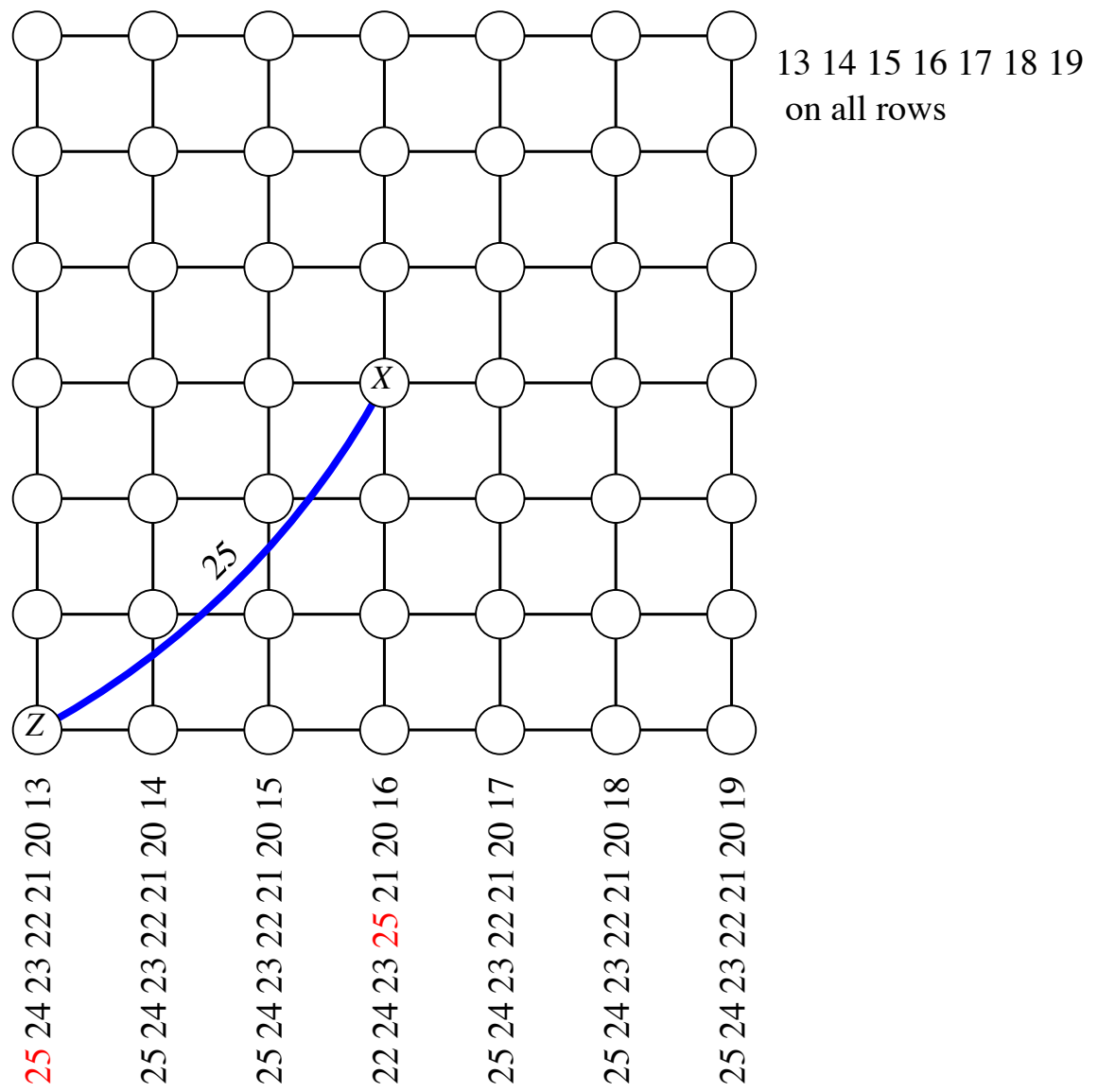

Figure 5. A class-2 coloring of $R_{7,7}$ with color 25 placed so it does not interfere with the bishop edge $X \leftrightarrow Z$; therefore 25 can replace 12 on that bishop edge, reducing the total color count to the desired 24.

We will build a class-2 coloring of $R_{n,n}$ using colors $2n-1,\ldots,4n-3$, which are unused in the bishop coloring (recall $\Delta(R_{n,n}) = 2n-2$). Color each row with $2n-1,\ldots,3n-2$ so that this order indicates the missing colors in each row (Fig. 5). Now color the columns by using the remaining colors together with the appropriate missing color; e.g., the first column gets colors $3n-1,\ldots,4n-3$ together with $2n-1$. Arrange the column coloring so that the missing colors at the rows are in reverse numerical order (as in Fig. 5), except that, in the central column, color $4n-3$ is missing at the central vertex $X$. Then we can finish by replacing color $2n-2$ on $X \leftrightarrow Z$ in the bishop coloring by color $4n-3$ (25 in Fig. 5). So the total number of colors used is now $4n-2$, and combining the two colorings gives a class-1 coloring of $Q_{n,n}$. In fact, it is a also an extremal coloring (see end of Section 2), as the rarest color appears only once. □

Returning to the general edge-coloring questions left open by Theorem 6, the most natural conjecture is that $Q_{m,n}$ is class 1 whenever it is not overfull. The first cases are: $Q_{3,n}$, $3 \le n \le 11$; $Q_{5,n}$, $5 \le n \le 69$; $Q_{7,n}$, $7 \le n \le 207$; $Q_{9,n}$, $9 \le n \le 457$; and $Q_{11,n}$, $11 \le n \le 851$.

**Conjecture 2.** The queen graph $Q_{m,n}$ is class 2 iff it is overfull.

We have some positive steps toward Conjecture 2. A first step was a generalization of the $Q_{m,m}$ case that combined bishop and rook colorings and worked for $m \le n \le 2m-1$; we omit the details because Theorem 7 in Section 4 uses more delicate arguments to get the much stronger result that $Q_{m,n}$ is class 1 when $m \le n \le \frac{1}{2}\left(m^2 - 3\,m + 2\right)$. For small values (e.g, $m = 3, 5$) the quadratic result is not as good as the $2m - 1$ result, but that is not a problem because various computations, which we describe in a moment, show that Conjecture 2 is true for $m \le 9$.

If $f(m) = \frac{1}{3}\left(2m^3 - 11m + 18\right)$, then $Q_{m,f(m)}$ is "just overfull" [SSTF, p.71], in that $n_e = \Delta\left\lfloor \frac{1}{2} n_v \right\rfloor + 1$. The *Just Overfull Conjecture* [SSTF, p. 71] states that for any simple graph $G$ such that $\Delta(G) \ge \frac{1}{2} n_v(G)$, $G$ is just overfull iff $G$ is "edge-critical" (meaning, $\chi'$ decreases upon the deletion of any edge). Computations show that $Q_{3,13}$ is edge-critical. Since deletion of a single queen edge cannot reduce the maximum degree, this is the same as saying that the deletion of any queen edge leads to a class-1 graph. So we have the following additional conjecture about the structure of queen edge colorings.

**Conjecture 3.** The queen graph $Q_{m,n}$ is just overfull iff it is edge critical.

An algorithm based on Kempe-style color switches yielded class-1 colorings for queen graphs that verify Conjecture 2 for $m = 3, 5, 7, 9$, and for $m = 11$ up to $n = 551$. A straightforward bootstrapping approach, where $Q_{m,n}$ was used to generate a precoloring for $Q_{m,n+2}$ and random Kempe color-switches were used to resolve impasses, worked for $m = 3$ and 5 (and 7 up to $Q_{7,199}$); an example of a class-1 coloring of $Q_{3,7}$ is in Figure 6; it was found by a method similar to the general Kempe method, but with an effort to find an extremal bishop coloring, which is shown at top with white the rarest color. But a more subtle method yielded a much faster algorithm, which resolved Conjecture 2 for $m = 7$ and 9 (and 11 up to $Q_{11,559}$, and also $Q_{11,707}$). The largest case required the coloring of 2,861,496 edges! This faster algorithm uses an explicit method to get a $\Delta + 1$ coloring (e.g., one can combine optimal colorings of the rook and bishop subgraphs) and then Kempe-type switches to eliminate the least popular color. This last step is based on a local search method that assigns a heuristic score to the possible switches and chooses the one with the highest score. This approach is quite general, using no information specific to the queen graph (except the speedy generation of the initial $\Delta + 1$ coloring, a task that can also be done via the algorithm inherent in Vizing's proof that a coloring in $\Delta + 1$ colors always exists).

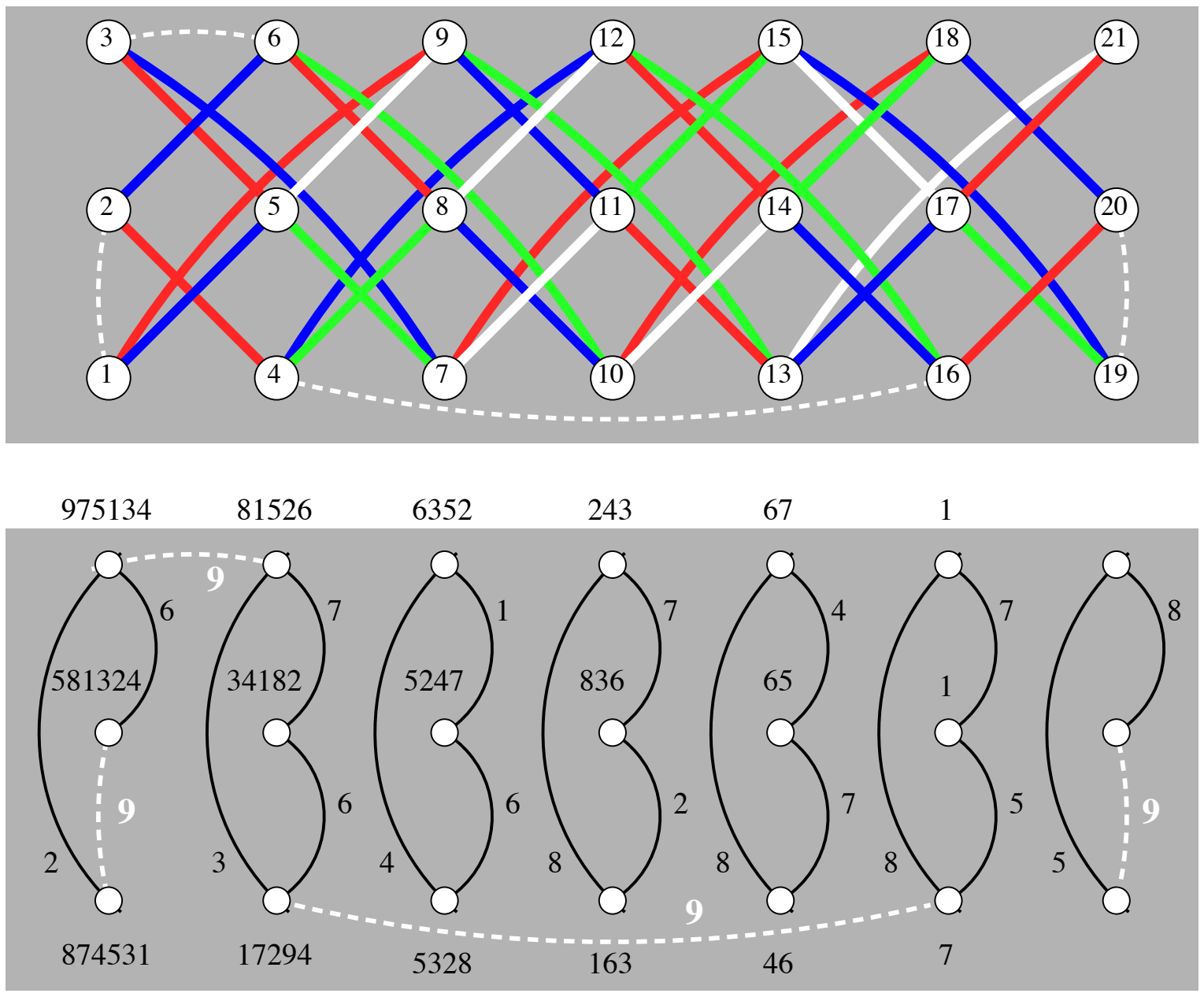

Figure 6. Top: An extremal edge coloring of $B_{3,7}$; there are four colors and white is avoided at vertices 1, 2, 3, 4, 6, 16, 18, 19, and 20. The dashed arcs indicated how white can then be used on four rook edges. Take the four colors, starting with white, to be 9, 10, 11, 12. Bottom: A class-1 coloring using 1 through 8 of the rook graph $R_{3,7}$ less the four edges from (a); only the vertical edges are shown as arcs in the edge-deleted graph. The four dashed white edges get the shared color, 9. The three sets of horizontal labels indicate edge colors on the horizontal edges, moving to the right. This rook coloring combines with the bishop coloring to yield a class-1 coloring of $Q_{3,7}$ using 1 through 12.

# 4. Quadratic Class 1 Queen Graphs

In this section we show how a certain multigraph defined from the canonical bishop coloring can help prove that many queen graphs are class 1; the method can be called the *ladder-and-multicycle method*. Throughout this section $m$ and $n$ are odd, $m \le n$, and $k = \frac{1}{2}(m-1)$. The main result, proved in Section 4.2, is the following.

**Theorem 7.** $Q_{m,n}$ is class 1 for all $m \le n \le \frac{1}{2}\left(m^2 - 3\,m + 2\right)$.

## 4.1 The Derived Multicycle of a Bishop Coloring

We will here use only the canonical path-based class-1 coloring of the bishop graph $B_{m,n}$, as described in Theorem 3. From such a coloring, we can define a derived multigraph; edge-coloring information about the multigraph can yield edge-coloring information about the corresponding queen graph $Q_{m,n}$. The multigraph is in fact a *multicycle*, by which we mean a multigraph on vertex-set $V$ with edges being edges of the associated simple cycle on $V$. Recall from Theorem 3 that the last two colors ($2m-3$ and $2m-2$; in this section we use cyan for color $2\,m-2$) in the canonical bishop coloring are used on the subgraph $G_k^-$, which consists of paths that start with edges of positive slope and length $k$, then negative slope edges of length $k+1$, then positive slope edges of length $k$, and so on (Fig. 2, bottom right). We call such a path a *special path*.

**Definition 1.** For any bishop graph $B_{m,n}$, let $\hat{B}_{m,n}$, the *derived multicycle*, be the multigraph on vertices $\{1, 2, \ldots, m\}$ given in the order $\{1 + j\,k : j = 0, \ldots, m-1\}$ (where the numbers are reduced mod $m$ starting from 1); the edges arise from the $(2m-2)$-colored bishop edges: the bishop edge $(x_1, y_1) \bullet\!\!-\!\!\bullet (x_2, y_2)$ induces the edge $y_1 \bullet\!\!-\!\!\bullet y_2$ in $\hat{B}_{m,n}$. The case of $\hat{B}_{5,11}$ is shown in Figure 7.

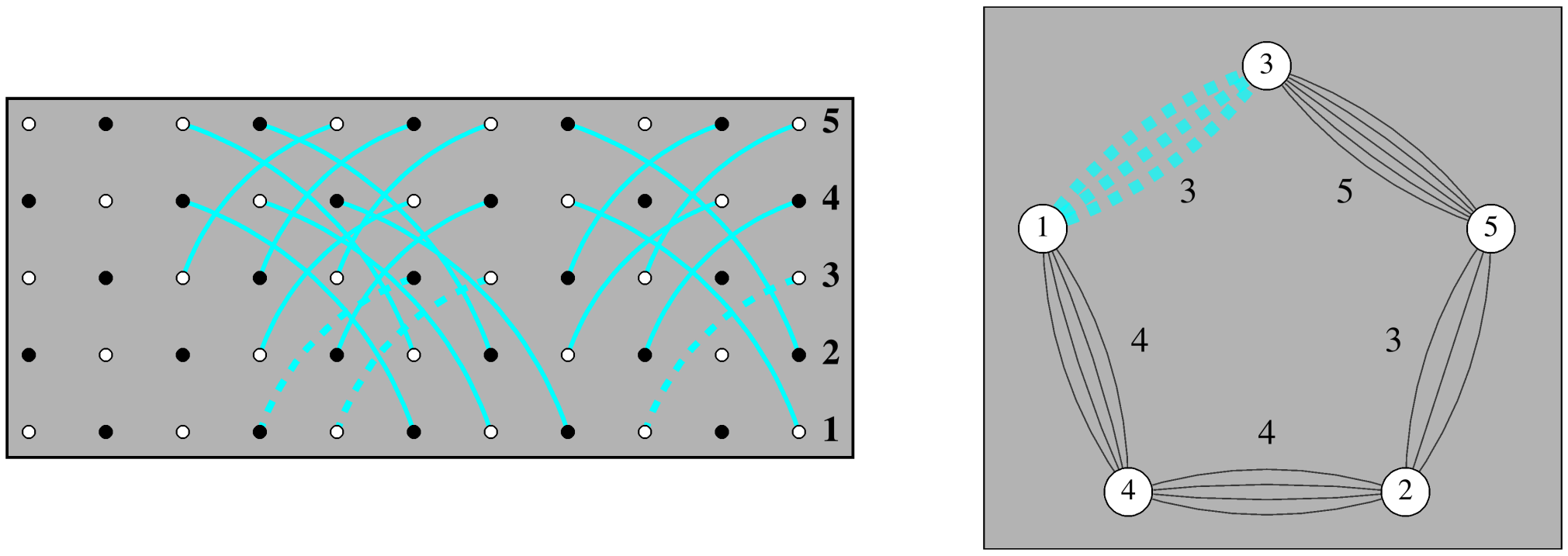

Figure 7. The last color (color $8 = 2\cdot 5 - 2$) of the canonical coloring of $B_{5,11}$ (left) and the derived multicycle $\hat{B}_{5,11}$ with edge-multiplicities shown. The dashed edges correspond.

The multiplicities of the edges in $\hat{B}_{m,n}$ ((3, 5, 3, 4, 4) in Fig. 7) play a key role in the proof that follows. The critical parameters of $\hat{B}_{m,n}$ are $\Delta$ and $\chi'$, and the minimum multiplicity $\mu^-$. We also use the maximum multiplicity $\mu^+$ and $\sigma_{m,n}$, the edge-count of $\hat{B}_{m,n}$ (i.e., the number of cyan edges in the canonical coloring). Now here is the key result that relates the chromatic index of $\hat{B}_{m,n}$ to that of $Q_{m,n}$.

**Proposition 8 (The Ladder-and-Multicycle Method).** Suppose $\hat{B}_{m,n}$ can be edge-colored using $n-1$ colors. Then $Q_{m,n}$ is class 1. That is, $\chi'(\hat{B}_{m,n}) \le n-1$ implies $\chi'(Q_{m,n}) = 3\,m + n - 4$.

**Proof.** We start with a special class-2 coloring (a "ladder coloring") of $R_{m,n}$ using colors $\{1, 2, \ldots, m-n-1\}$. Because the rook graph is regular, each vertex in any class-2 coloring misses exactly one color. We start by using 1 through $m$ on the columns (and some row edges), and will then use the $n-1$ colors $A = \{m+1, \ldots, m+n-1\}$ on the uncolored row edges. Each vertex in the leftmost column will use all colors in $A$, but the sequence of $n-1$ missing colors in each row excluding its leftmost vertex is a permutation of $A$; moreover, the colors in $A$ may be arranged so that, for each row, any preselected permutation is the missing-color permutation for that row. To define the ladder coloring, use colors 1 through $m$ on each column, ensuring that color $i$ is missing at vertices in the $i$th row. Now use 1 to color the horizontal edges in the bottom row that connect vertices in successive columns after the leftmost; i.e., the edges connecting the vertices in columns 2 and 3, columns 4 and 5, and so on. Do the same for row 2 but using color 2, and so on (Fig. 8). We have now used $m$ colors to color all vertical edges and the edges of one maximum matching in each row. But each row is a $K_n$, and $K_n$ minus any maximum matching can be colored with $n-1$ colors (Prop. 1). Thus the colors in $A$ suffice to color all uncolored horizontal edges. The vertices in the leftmost column see all the colors in $A$, while the remaining vertices (which already have edges colored 1 through $m$) each miss exactly one color in $A$. Thus the missing colors in each left-deleted row form a permutation of $A$, and it is clear from the construction that the $A$-colors can be arranged independently in the rows, so that any set of $m$ permutations can be assumed to be the missing-color permutations on the rows, excluding the leftmost vertices.

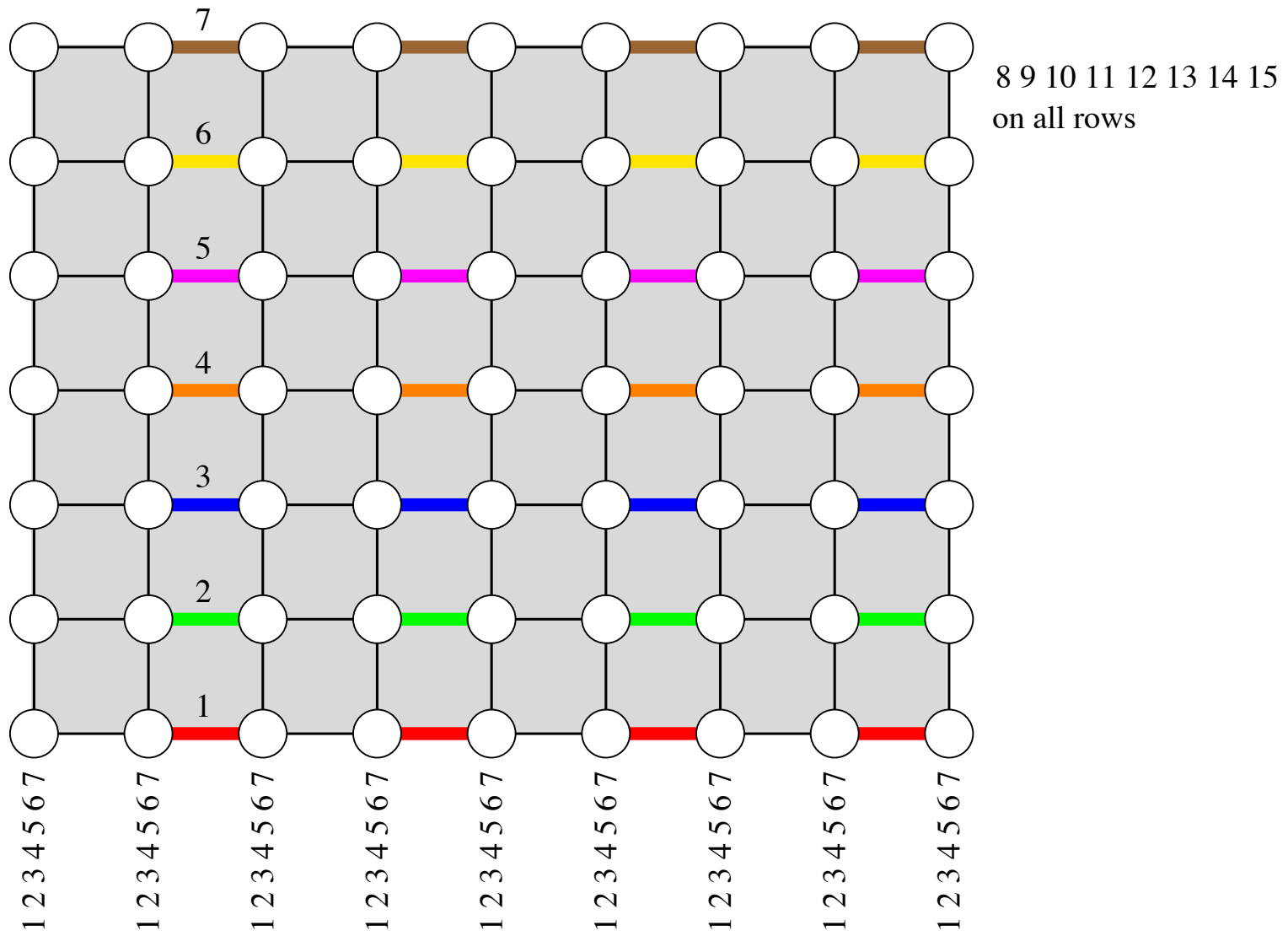

Figure 8. A class-2 ladder coloring of $R_{7,9}$ using 15 colors. The column colors are listed in missing order.

The proof of the theorem now proceeds as follows. We assume that the bishop graph is colored in the canonical way (Thm. 3) using colors from 1 to $2m-2$. Let $\xi$ be the hypothesized edge-coloring of $\hat{B}_{m,n}$ using colors $\{m+1, \ldots, m+r\}$, $r \le n-1$; this is a subset of $A$. If $\xi$ assigns $k$ to edge $e = y \multimap y'$ of $\hat{B}_{m,n}$, assign $k$ as a missing color to the endpoints of $e$ viewed as a bishop edge; that is, if $e$ arises from the $(2m-2)$-colored bishop edge $(x, y) \multimap (x', y')$, assign $k$ to be used as a missing rook color at both vertices $(x, y)$ and $(x', y')$. Because $\xi$ is a proper multigraph coloring of $\hat{B}_{m,n}$, no $k$ will be assigned to more than one vertex in any row; and because no bishop edge incident with the leftmost column gets color $2m-2$, no missing color will be assigned to vertices in the leftmost column. We therefore get, for each row, an injection from $\{m+1, \ldots, m+r\}$ to the vertices of that row excluding its leftmost vertex, and these maps can be extended to full permutations of $A$ arbitrarily. Since, in the class-2 rook coloring, we can arrange the missing $A$-colors to match these missing-color permutations, we can now recolor each $(2m-2)$-colored bishop edge with the common missing rook color at its endpoints, thus eliminating $2m-2$ as a bishop color. So now the two colorings combine to give a coloring of $Q_{m,n}$ with color count equal to $\Delta(B_{m,n}) - 1 + \Delta(R_{m,n}) + 1$, which is $\Delta(Q_{m,n})$. □

We will use Proposition 8 to obtain an infinite family of queen class-1 colorings, but first we need careful analysis of multicycles in general and also of the particular multicycles $\hat{B}_{m,n}$.

Next we turn to an analysis of the chromatic index of multicycles. There might be an efficient formula or algorithm that can go quickly from the multiplicity vector for a multicycle to its chromatic index, but there are several types of things that can happen and we do not have such a general method. It seems unlikely, but perhaps the problem is $\mathcal{NP}$-hard! But for multicycles with certain nice structure it is easy to get an exact formula or a tight upper bound. Recall the two classic, fully general bounds for multigraphs: Vizing's bound is $\chi'(G) \le \Delta(G) + \mu^+(G)$ and Shannon's bound is $\chi'(G) \le \left\lfloor \frac{3}{2}\Delta(G) \right\rfloor$. Let $C_{m,a}$ denote the regular multicycle of length $m$ with $a$ edges in each group. Determining the chromatic index for regular multicycles is easy.

**Proposition 9 (Regular Multicycle Chromatic Index).** $\chi'(C_{m,a}) = 2\,a + \left\lceil \frac{2\,a}{m-1} \right\rceil$.

**Proof.** View $C_{m,a}$ as a collection of $a$ simple cycles and partition them into groups of size $k = (m-1)/2$, or less for the last group. Each group can be edge-colored using two colors for each cycle plus one extra color that is used on all the cycles in the group, spreading the extra color around a maximum matching in the $m$-cycle so that the edges with this extra color are disjoint (Fig. 8). The total color count is $2a + \left\lceil \frac{a}{k} \right\rceil$. This upper bound is sharp because $\rho(C_{m,a}) = k$. If the color count was less than the upper bound, the edge count would be at most $k\left(\left\lceil \frac{a}{k} \right\rceil + 2\,a - 1\right)$; using the identity $k \left\lceil \frac{a}{k} \right\rceil \leq a + k - 1$ simplifies this to $m\,a - 1$, one less than the number of edges. □

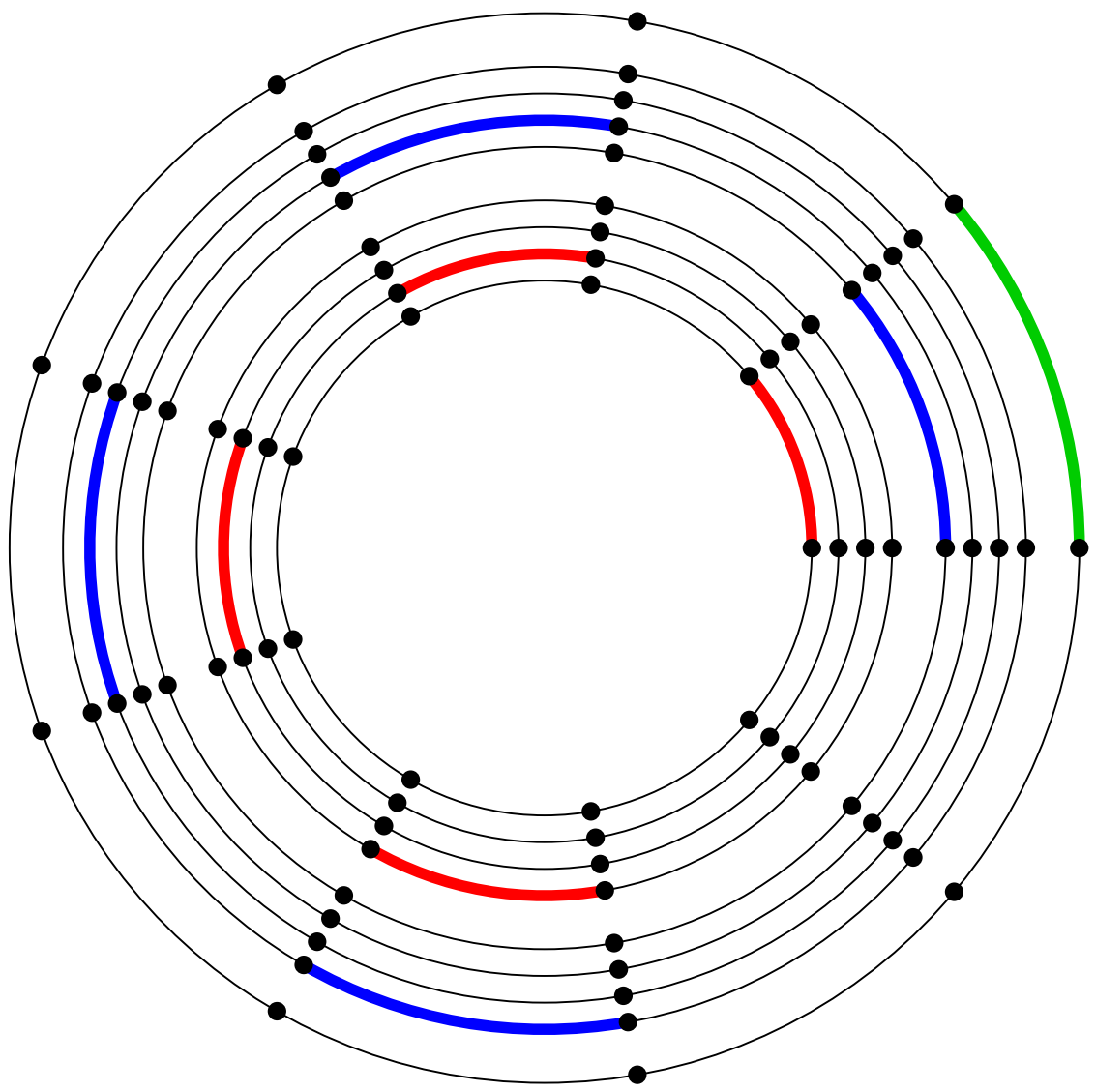

Figure 8. The regular multicycle $C_{9,9}$ (shown split apart into 9 cycles) can be edge-colored using 21 colors. Each cycle gets 2 colors (for 18), with three shared colors (one for each group of four or less; red, green, blue) each placed in a matching.

If $m$ is even, then $\chi'(m, a) = 2a$, but this is irrelevant to our work. More important here is the simple case is a *multipath*: a multicycle that has at least one 0 multiplicity.

**Lemma 10.** If $G$ is a multipath, then $\chi'(G) = \Delta(G)$.

**Proof.** Enumerate the edges in the order they appear in the path as $\{e_i\}$ and assign color $i \pmod{\Delta}$ to $e_i$. □

The preceding cases lead to a simple upper bound for any multicycle.

**Corollary 11 (Multicycle Chromatic Index Bound).** Let $G$ be any multicycle with vertex count $m$ and multiplicity vector $\vec{M}$. Then $\chi'(G) \leq \chi'(C_{m,\mu^-}) + \Delta(G) - 2\mu^-$. Using $k$ for $\frac{1}{2}(m-1)$ as usual, this becomes
$\chi'(G) \leq \left\lceil \frac{\mu^-}{k} \right\rceil + 2\mu^- + \Delta(G) - 2\,\mu^- = \Delta(G) + \left\lceil \frac{\mu^-}{k} \right\rceil$.

**Proof.** Split $G$ into the regular "kernel" — the multicycle $C_{m,\mu^-}$ — and the "residual", which is a *multipath* (because kernel removal leaves a 0 multiplicity) and so has chromatic index $\Delta(G) - 2\mu^-$ by Lemma 10. The kernel is colorable as in the Regular Multicycle Coloring Theorem, and summing the two yields the claimed bound. □

Compare the preceding bound to the general Vizing bound: $\chi' \leq \Delta + \mu^+$. For multicycles, $\chi' \leq \Delta + \left\lfloor \frac{\mu^-}{k} \right\rfloor$. The values $\mu^\pm(\hat{B}_{m,n})$ do not differ by much (in general of course they can differ arbitrarily), but division by $k$ is a big improvement.

Recall the trivial lower bound $\left\lceil \frac{n_e(G)}{\rho(G)} \right\rceil \le \chi'(G)$. For our bishop multicycles this becomes $\left\lceil \frac{\sigma(\hat{B}_{m,n})}{k} \right\rceil \le \chi'(\hat{B}_{m,n})$. In fact, it appears that this lower bound is always the exact value of the chromatic index for the derived multigraph. We have checked this for $m \le 19$ and $n \le 199$. To find the true chromatic index in these cases, we used a composite algorithm that we now describe.

First one can try a greedy enumeration procedure. Set $d = \Delta(G)$ and simply enumerate the edges in the order they appear around the multicycle as $\{e_i\}$; assign color $i \pmod d$ to $e_i$. If this process never leads to a conflict, $d$ is an upper bound on $\chi'$. If a conflict does arise, increment $d$ by 1 and start over; continue until a proper coloring is found. If the resulting upper bound agrees with $\tau$, the trivial lower bound, $\chi'$ is proved to be $\tau$.

Then there is the upper bound of Corollary 11; call it $\tau_1$. If $\tau_1 = \tau$, then the chromatic index is $\tau$. One can sometimes reduce $\tau$ by 1 to get an *improved upper bound*. For example, if $\mu \equiv 1 \pmod k$ then the last cycle of the kernel is in a group of 1 and requires 3 colors; with the third color available to be placed arbitrarily. But if the number of major vertices in the residual is 1 or the number of isolated edges is 1, then that third color can be shared with the residual. Also if $\mu \not\equiv 0 \pmod k$ and the residual has only one edge, or consists only of isolated edges, then again a reduction by 1 is allowable.

If the preceding two bounds fail, split $\hat{B}_{m,n}$ into the kernel and residual and use the preceding improved upper bound method on the residual and the formula of Proposition 9 for the kernel. If they sum to $\tau$, again we have succeeded in finding $\chi'$. These methods are sufficient to find the chromatic index in the 756 cases we have examined, and so it is reasonable to conjecture that the chromatic index of $\hat{B}_{m,n}$ agrees exactly with the lower bound $\left\lceil \frac{2\sigma}{m-1} \right\rceil$ in all cases. This lower bound is not exact for general multicycles; for multiplicities (0, 0, 0, 1, 2) the lower bound is 2 but $\chi' = 3$.

**Conjecture 4.** $\chi'(\hat{B}_{m,n}) = \left\lceil \frac{2}{m-1} n_e(\hat{B}_{m,n}) \right\rceil$.

There are many many patterns in the data one can compute for the derived multicycle of $B_{m,n}$; key parameters are the total edge count $\sigma$, the minimum multiplicity $\mu^-$, and the maximum degree $\Delta$. The next conjecture summarizes the results of many computations. Figure 8 presents some evidence for Conjecture 5, and also shows the periodicity in the edge counts of $\hat{B}_{m,n}$ that appears to arise in all cases.

**Conjecture 5.** For odd $m, n$, with $m \le n$, $\frac{1}{2} m\,n - \left(\frac{1}{2} m^2 - 1\right) \le \sigma_{m,n} \le \frac{1}{2} m\,n - \frac{1}{4}\left(m^2 + 1\right)$.

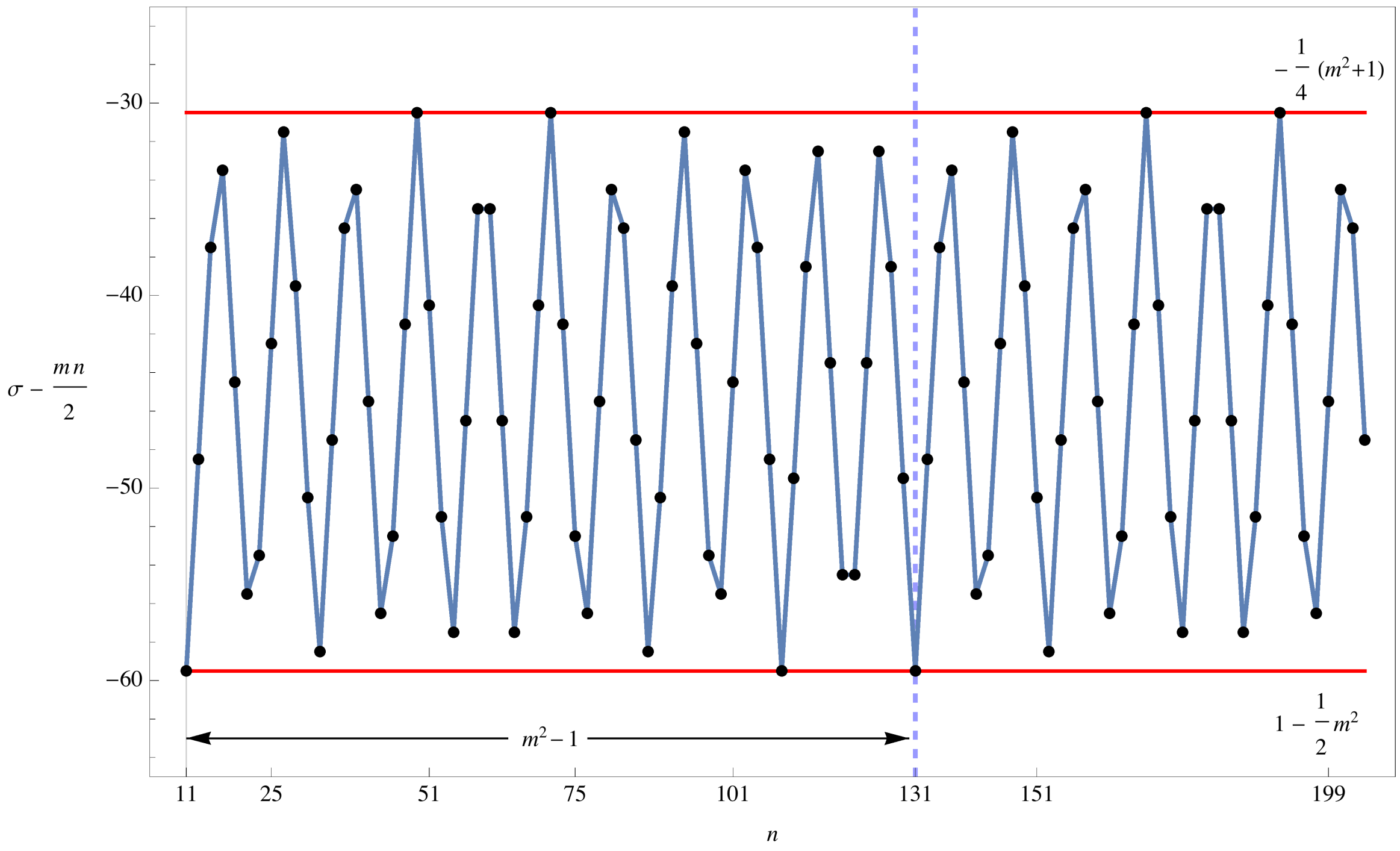

Figure 8. Computed values of $\sigma_{m,n} - \frac{m\,n}{2}$, where $m = 11$ and with the red lines being $-\frac{1}{4}\left(m^2 + 1\right)$ (upper) and $1 - \frac{1}{2}m^2$ (lower). This illustrates the bounds of Conjecture 5. Also note the periodicity of this reduced data set with period $m^2 - 1$ (dashed blue line).

## Limits of the Ladder-and-Multicycle Method

Although we will use the ladder-and-multicycle method to prove Theorem 7, we note here that it cannot be used to prove the full Conjecture 2. Suppose $2m - 1 \le n$; then the number of major vertices in $B_{m,n}$ is $m\,n - (3m + 1)\,k$ (easy exercise: $n - 2\,m\,(m - 1) + 2\,k + 4\,\Sigma_{i=1}^{k-1}\,i$); therefore, for such $n$, any bishop color class from a class-1 coloring must have size at least $a = \frac{1}{2}\,(m\,n - (3\,m + 1)\,k\, - 1) + 1$. Now suppose that $\chi'(\hat{B}_{m,n}) \le n - 1$; then $\hat{B}_{m,n}$ is covered by at most $n - 1$ matchings, each of size at most $k$, and thus $\hat{B}_{m,n}$ can have at most $b = k(n - 1)$ edges. But this is the same as the size of the $(2m - 2)$-color class, so we must have $a \le b$, which simplifies to $n \le \frac{3}{2}\,m^2 - 2m - \frac{1}{2}$. This shows that the ladder-and-multicycle method can work only up to this quadratic function. In fact, assuming certain observed patterns hold, it appears that the method can work only to about $m^2$. For assume Conjecture 5; then $\sigma_{m,n} \ge \frac{1}{2}\,m\,n + 1 - \frac{1}{2}\,m^2$. Therefore $\chi'(\hat{B}_{m,n}) \ge \left(\frac{1}{2}\,m\,n + 1 - \frac{1}{2}\,m^2\right)/k$. Now for this to be at most $n - 1$ means $n \le m^2 - m - 1$, so this is a likely bound for the ladder-and-multicycle. In the next section we will show that the method can be proved to work for values of $n$ near $\frac{1}{2}\,m^2$.

## 4.2 The Fine Structure of the Canonical Bishop Coloring

We can use the color-sharing theorem and a detailed study of the properties of the last color in the canonical bishop coloring to prove the queen coloring result of Theorem 7.

**Proof of Theorem 7.** Since we have already shown that $Q_{m,n}$ is class 1 if either $m$ or $n$ is even, we assume $m$ and $n$ are odd and write $m = 2\,k + 1$. Recall (Thm. 3) that in $B_{m,n}$ only two colors (one of them being cyan) were needed to color all special paths. As we are at liberty to start each special path with either color, we assume that no special path starts with cyan.

Let $C$ be the set of cyan edges in all special paths. Because $k$ and $m$ are relatively prime, the derived multigraph is a multicycle. The following result will yield Theorem 7 as a consequence of the Ladder-And-Multicycle Method (Prop. 8) and the Multicycle Chromatic Index Bound. It is not hard to see that $n - 2\,k \le \Delta(\hat{B}_{m,n})$. Proposition 12, the final step in the proof, puts an upper bound on the maximum degree; so we see that, as $m$ is fixed and $n$ rises, the derived multicycle is close to being regular.

**Proposition 12.** When $m$ and $n$ are odd, $m \le n$, $m = 2\,k + 1$, we have $\Delta(\hat{B}_{m,n}) \le n - k$.

We can then complete the proof of Theorem 7 as follows. For any multicycle $G$, $\mu^{-}(G) \le \frac{1}{2}\,\Delta(G)$, so by Corollary 11 we have $\chi'(\hat{B}_{m,n}) \le n - k + \left\lceil \frac{n-k}{2\,k} \right\rceil$. By assumption, $n \le \frac{1}{2}\left(m^2 - 3m + 2\right) = k(2k - 1)$, and therefore $\frac{n-k}{2\,k} \le k - 1$ so that $\left\lceil \frac{n-k}{2\,k} \right\rceil \le k - 1$ and $\chi'(\hat{B}_{m,n}) \le n - 1$. Theorem 8 now concludes the proof.

We now prove Proposition 12, starting with some definitions: call the leftmost vertex of a special path in a bishop graph an *initial vertex* and let $\ell(i, j)$ be the length (i.e., edge count) of the special path that starts at the initial vertex $(i, j)$. Furthermore, call an initial vertex $(i, j)$ for which $\ell(i, j)$ is even an *even vertex* and one for which $\ell(i, j)$ is odd an *odd vertex*. For bishop vertices $v$ viewed as points in the plane, we use $v \le w$ to mean that $v$ is first in the lexicographic ordering: $v_x \le w_x - 1$, or $v_x = w_x$ and $v_y \le w_y$. We break the proof into a series of claims.

**Claim 1.** Let $I(j)$ be the number of initial vertices in row $j$, and let $O(j)$ be the number of odd vertices in row $j$. Then the degree of row $j$ in $\hat{B}_{m,n}$ is $\deg(j) = n - I(j) - O(m + 1 - j)$.

**Proof.** Note that the degree of a vertex in $\hat{B}_{m,n}$ is the number of bishop vertices in the row represented by the multigraph vertex that are incident with an edge in $C$. Every vertex in $B_{m,n}$ that is not an endpoint of a special path lies on a cyan edge, so is counted toward the degree of the row in which it sits. Because we have arranged for all paths to start on the left without cyan, we eliminate all initial vertices in the row from the count; we also eliminate all vertices in the row that terminate an odd-length path, as they will not be incident with a cyan edge. Note that $\pi(i, j) = (n+1-i, m+1-j)$, a vertical reflection followed by a horizontal reflection, is an automorphism of $B_{m,n}$ that takes special paths to special paths (and odd length special paths to odd length special paths). Therefore if $v$ is an odd terminal point of path $p$, $\pi(v)$ is an odd vertex beginning path $\pi[p]$. Thus the number of vertices in row $j$ that terminate an odd-length path is the same as the number of odd (initial) vertices in row $m+1-j$, and the claim follows. □

**Claim 2.** $I(j) = k+1$ for $j \leq k$; $I(j) = k$ for $j > k$.

**Proof.** This follows from the fact that the set of initial points of all special paths is $\{(i, j) : (i, j) \leq (k+1, k)\}$; see Figure 9. □

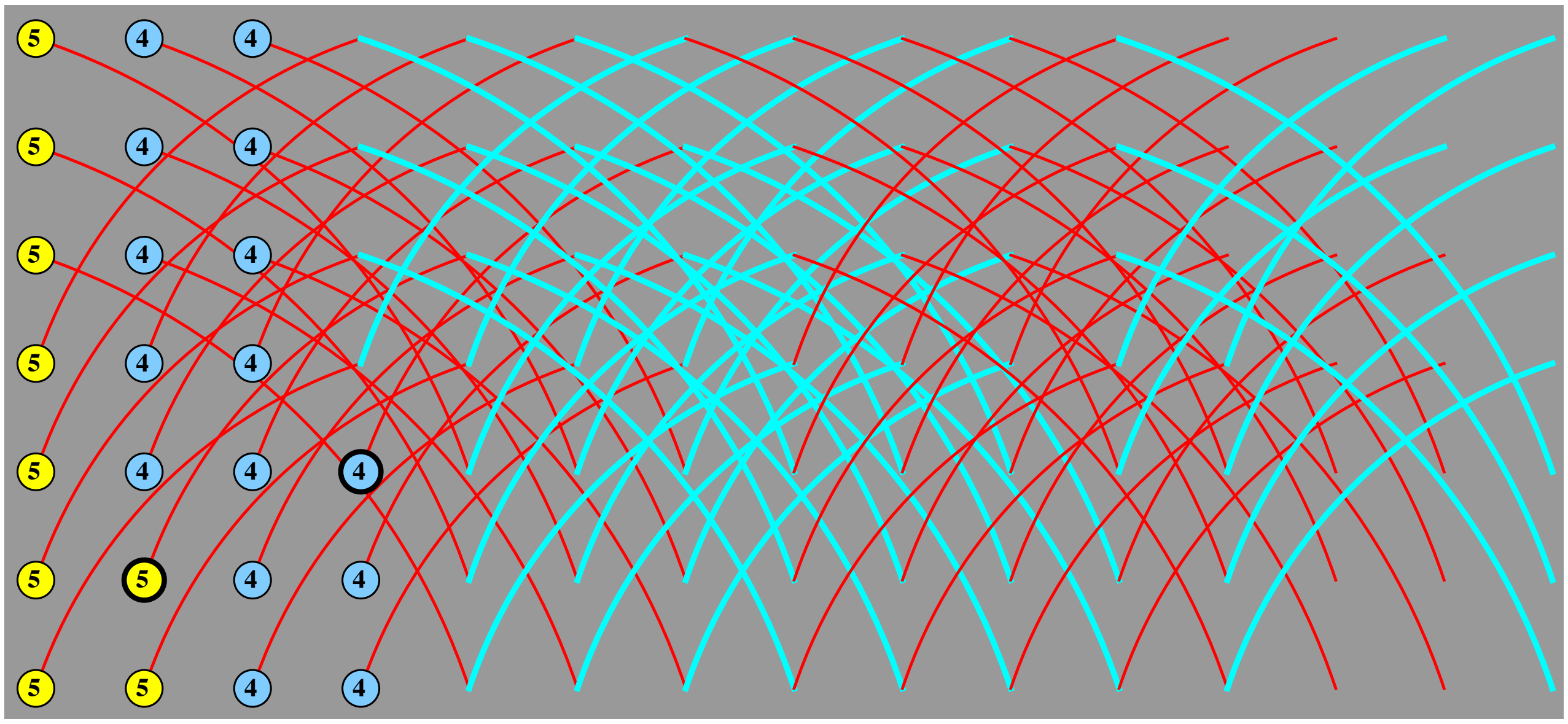

Figure 9. The values of $I(j)$ for $B_{7,15}$ are 4, 4, 4, 3, 3, 3 (Claim 2); these are all the vertices lexicographically below or equal to (4, 3). The yellow vertices are even, the blue ones odd.

**Claim 3.** There is at least one even vertex and at least one odd vertex.

**Proof.** Because $C$ is a matching, $e = v/2$ where $e = |C|$ and $v$ is the number of vertices in $B_{m,n}$ that are included in a $C$ edge. But the vertices included in such an edge are those that don't begin a special path or terminate an odd-length special path. There are $k(m+1)$ vertices that begin a special path, so letting $p$ be the number of odd vertices, we have $e = \frac{1}{2}(n\,m - k(m+1) - p)$. Because $n\,m$ is odd and $k(m+1)$ is even, this implies that $p$ must be odd, and thus there are an odd number of odd vertices, and hence also an odd number of even vertices. □

**Claim 4.** $\ell(i, j)$ is the largest integer $q$ such that $qk + i + \left\lfloor \frac{q\,k+j-1}{m} \right\rfloor \leq n$.

**Proof.** A special path moves $k$ units to the right (i.e., from $i$ to $i+k$) when it moves upwards and $k+1$ units to the right when it moves downwards, and only moves downwards if $c\,k + j \leq b\,m < (c+1)\,k + j$, where $c$ is the number of edges the path has moved along thus far, and $b\,m$ is any multiple of $m$. Thus a path that starts at $(i, j)$ and moves to the right $q$ edges ends at $(i', j')$, where $i' = q\,k + i + \lfloor (q\,k + j - 1)/m \rfloor$, and the claim follows.

**Claim 5.** There is an initial vertex $(i^*, j^*) < (k+1, k)$ such that $\ell(i, j) = \ell(1, 1)$ for all $(i, j) \leq (i^*, j^*)$, and $\ell(i, j) = \ell(1, 1) - 1$ for all initial vertices $(i, j) > (i^*, j^*)$.

**Proof.** See Figure 9 for an example where $(i^*, j^*) = (2, 2)$. Let $L = \ell(k+1, k)$. It follows easily from Claim 4 that $\ell(i, j)$ is decreasing with respect to lexicographic order; that is, $(i, j) \leq (i', j')$ implies $\ell(i, j) \geq \ell(i', j')$. Thus $L \leq \ell(1, 1)$, and the inequality must be strict by Claim 3. Also, by Claim 4, we have that $\ell(1, 1)$ is the largest integer $q$ such that $qk + \lfloor \frac{qk}{m} \rfloor \le n - 1$, while $L$ is the largest integer $q$ such that $q\,k + k + \lfloor \frac{q\,k+k-1}{m} \rfloor = (q+1)\,k + \lfloor \frac{(q+1)\,k-1}{m} \rfloor \le n - 1$. This implies that $\ell(1, 1) \le L + 1$. Thus $L \le \ell(1, 1) \le L + 1$, and the claim follows. More precisely, $(i^*, j^*)$ is the smaller of the largest even initial point and largest odd initial point, where comparisons are lexicographic. □

**Claim 6.** If $\ell(1, 1)$ is odd, $\Delta(\hat{B}_{m,n}) \le n - k$.

**Proof.** If $\ell(1, 1)$ is odd, Claim 5 implies that the odd vertices are all $(i, j) \le (i^*, j^*)$. Thus

$$O(j) = \begin{cases} i^* & \text{if } j \le j^* \\ i^* - 1 & \text{if } j > j^* \end{cases}$$

Then

$$O(m + 1 - j) = \begin{cases} i^* - 1 & \text{if } j < m + 1 - j^* \\ i^* & \text{if } j \ge m + 1 - j^* \end{cases}$$

By Claim 2, we therefore have

$$I(j) + O(m + 1 - j) = \begin{cases} i^* + k & \text{if } j \le \min(k + 1, m + 1 - j^*) \\ i^* + k - 1 & \text{if } k + 1 \le j < m + 1 - j^* \\ i^* + k + 1 & \text{if } m + 1 - j^* \le j < k + 1 \\ i^* + k & \text{if } j \ge \max(k + 1, m + 1 - j^*) \end{cases}$$

Thus $\Delta(\hat{B}_{m,n}) \le n - k - i^* + 1$ by Claim 1, and because $i^* \ge 1$, $n - k - i^* + 1 \le n - k$. □

We complete the proof of Propostion 12, and thus Theorem 7, with:

**Claim 7.** If $\ell(1, 1)$ is even, $\Delta(\hat{B}_{m,n}) \le n - k$.

**Proof.** If $\ell(1, 1)$ is even, Claim 5 implies that the odd vertices are all $(i, j)$ such that $(i^*, j^*) < (i, j) \le (k + 1, k)$. Thus

$$O(j) = \begin{cases} k + 1 - i^* & \text{if } j \le \min(j^*, k) \\ k - i^* & \text{if } k < j \le j^* \\ k + 2 - i^* & \text{if } j^* < j \le k \\ k + 1 - i^* & \text{if } j > \max(j^*, k) \end{cases}$$

Then

$$O(m + 1 - j) = \begin{cases} k + 1 - i^* & \text{if } j < m + 1 - \max(j^*, k) \\ k + 2 - i^* & \text{if } m + 1 - k = k + 2 \le j < m + 1 - j^* \\ k - i^* & \text{if } m + 1 - j^* \le j < k + 2 \\ k + 1 - i^* & \text{if } j \ge m + 1 - \min(j^*, k) \end{cases}$$

By Claim 2, we therefore have

$$I(j) + O(m + 1 - j) = \begin{cases} m - i^* + 1 & \text{if } j < m + 1 - \max(j^*, k) \text{ and } j \le k \\ m - i^* & \text{if } j < m + 1 - \max(j^*, k) \text{ and } j > k \\ m - i^* + 1 & \text{if } k + 2 \le j < m + 1 - j^* \\ m - i^* & \text{if } m + 1 - j^* \le j \le k \\ m - i^* - 1 & \text{if } m + 1 - j^* \le j = k + 1 \\ m - i^* & \text{if } j \ge m + 1 - \min(j^*, k) \end{cases}$$

Note that if $i^* = k + 1$, $j^*$ must be less than $k$ (by Claim 3), and thus the fifth of these cases $(m + 1 - j^* \le k + 1)$ cannot occur. In that case, $\Delta(\hat{B}_{m,n}) \le n - m + i^* = n - m + k + 1 = n - k$ by Claim 1. If $i^* < k + 1$, we similarly have $\Delta(\hat{B}_{m,n}) \le n - m + i^* + 1 \le n - k$, and the claim follows. □

The quadratic lower bound, $n \ge \frac{1}{2}\left(m^2 - 3\,m + 4\right)$, for class-2 queen graphs of Theorem 7 provides substantial support for the queen chromatic index conjecture (Conjecture 2). With more work, it is likely that the proof can be extended to yield a slightly higher bound, namely $\frac{1}{2}(m^2 + 2m - 3)$, but as shown at the end of Section 4.1, the ladder-and-multicycle method breaks down for $n > \frac{1}{2}(3m^2 - 4m - 1)$, and thus a quadratic bound is the best we can achieve in this manner. Therefore a proof that $Q_{m,n}$ is class 1 all the way up to the cubic bound of Conjecture 2 would seem to require an altogether different approach of even greater intricacy. This highlights the surprising difficulty of the Vizing classification problem for queen graphs. It is well known that the general classification problem is hard ($\mathcal{NP}$-complete), but one would expect that queen graphs — being a union of rook and bishop graphs whose classifications are relatively straightforward — would be amenable to a complete classification. So as we leave behind the quadratic bound obtained here, it is hard to resist the thought that we may be entering terrain of intractable complexity. If so, the computational evidence supporting Conjecture 2 seems all the more remarkable and may be the best we can hope for.

## 5. Mycielski Graphs

The *Mycielskian* $\mu(G)$ of a graph $G$ with vertex set $X$ is an extension of $G$ to the vertex set $X \cup Y \cup \{z\}$, where $|Y| = |X|$ and with new edges $z \leftrightarrow y_i$ for all $i$ and $x_i \leftrightarrow y_j$ for each edge $x_i \leftrightarrow x_j$ in $G$ (see Fig. 10). The *Mycielski graphs* $M_n$ are formed by iterating $\mu$ on the singleton graph $M_1$, but ignoring the isolated vertex that arises in $\mu(M_1)$. Thus $M_1 = K_1$, $M_2 = K_2$, $M_3 = C_5$, and $M_4$ is the Grötzsch graph (Fig. 11). They are of interest because $M_n$ is a triangle-free graph of chromatic number $n$ having the smallest possible vertex count. Fisher et al [FMB] proved that if $G$ is Hamiltonian, then so is $\mu(G)$. We extend that to a Hamilton-connected (HC) result, provided $n_v(G)$ is odd. This is sufficient to show that all Mycielski graphs $M_n$, except $M_3 = C_5$, are HC.

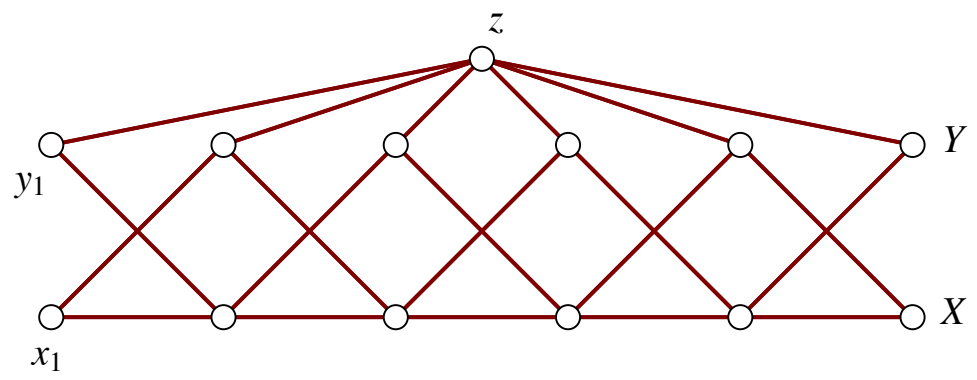

Figure 10. The Mycielskian of a 6-path.

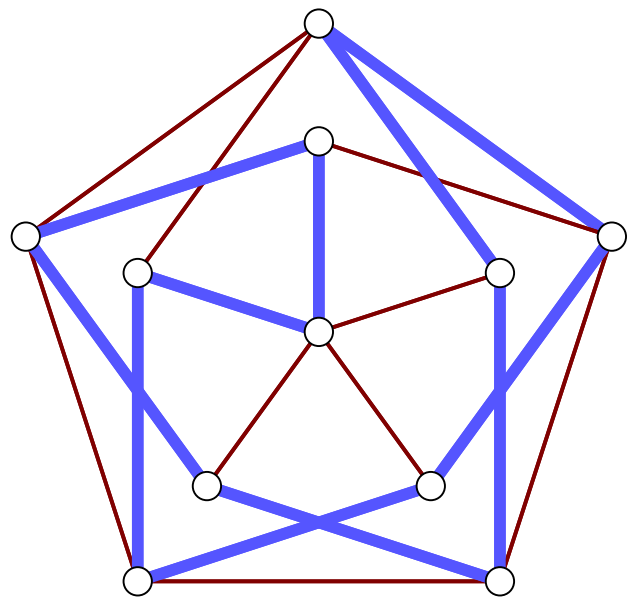

Figure 11. A Hamiltonian cycle in the Mycielski graph $M_4$, which is the Grötzsch graph.

The key result is the following.

**Theorem 13.** If $G$ is an odd cycle, then $\mu(G)$ is Hamilton-connected.

We will prove this shortly; note that it yields the fact that $\mu$ preserves HC for graphs with an odd number of vertices.

**Corollary 14.** If $G$ is Hamilton-connected and $n_v(G)$ is odd, then $\mu(G)$ is Hamilton-connected.

**Proof.** We may skip the trivial case that $G$ has one vertex. Therefore $G$ is Hamiltonian, with Hamiltonian cycle $C$. Since $\mu(G)$ contains $\mu(C)$ as an edge subgraph, and since $\mu(G)$ is HC by Theorem 13, so is $\mu(G)$.

Theorem 13 does not extend to the even case.

**Proposition 15.** If $G$ is an even cycle, then $\mu(G)$ is not HC.

**Proof.** It is easy to use parity to show that there is no Hamiltonian path from any vertex in $X$ to $z$. This is because such a path can get to $z$ only from $Y$, and hence must alternate from $X$ to $Y$; but then if the path starts at $x_1$ it can never visit $x_i$, where $i$ is even. □

Even cycles are not HC, so the negative result for even cycles does not mean that HC-preservation fails in general for even graphs (an exception being $K_2$: $\mu(K_2)$ is a 5-cycle, which is not HC, even though $K_2$ is HC). Computations support the following strengthening of Corollary 14, but some new ideas are needed.

**Conjecture 6.** If $G$ is Hamilton-connected and not $K_2$, then $\mu(G)$ is Hamilton-connected.

Easy computation shows that $M_4$, the 11-vertex Grötzsch graph, is HC and so Corollary 14 means that all Mycielski graphs are HC, except the 5-cycle $M_3$.

**Corollary 16.** The Mycielski graph $M_n$ is Hamilton-connected iff $n \neq 3$.

The stronger assertion that $\mu(G)$ is HC whenever $G$ is Hamiltonian is false. Counterexamples include $C_4$, $K_{3,3}$, $K_{1,1,2}$, $\text{Grid}_{2,3}$. Indeed, the first two here are Hamilton-laceable, but their Mycielskians are not HC.

**Proof of Theorem 13.** Assume that the cycle $G$ has $n$ vertices, given in cyclic order as $X = \{x_i\}$, where $n$ is odd. Then $\mu(G)$ has as vertices $X$, and also $Y = \{y_i\}$ and a single vertex $z$. The subgraph corresponding to $Y \bigcup \{z\}$ forms a $K_{1,n}$. Then, as in [FMB], we have the following Hamiltonian cycle in $\mu(G)$ (see Fig. 12):

$$C = y_1 \leftrightarrow x_2 \leftrightarrow y_3 \leftrightarrow x_4 \leftrightarrow \cdots \leftrightarrow x_{n-1} \leftrightarrow y_n \leftrightarrow x_1 \leftrightarrow x_n \leftrightarrow y_{n-1} \leftrightarrow \cdots \leftrightarrow y_4 \leftrightarrow x_3 \leftrightarrow y_2 \leftrightarrow z \leftrightarrow y_1$$

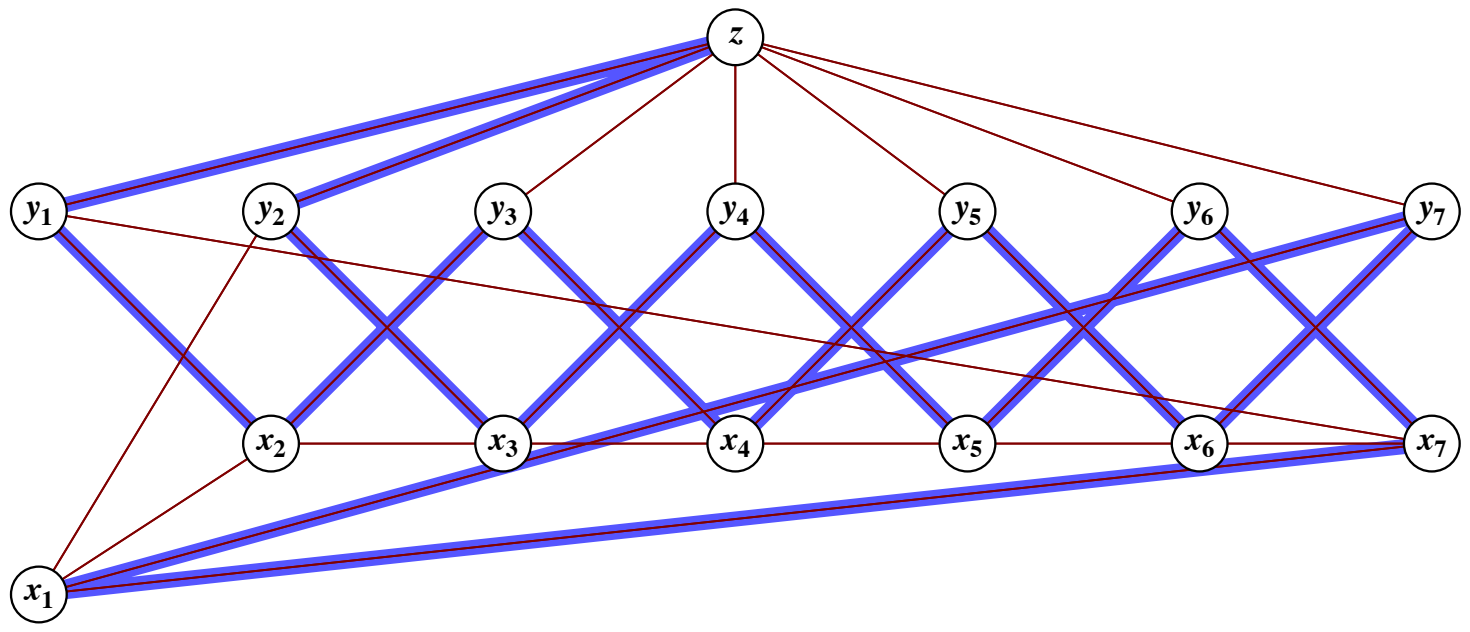

Figure 12. A Hamiltonian cycle in the Mycielskian of an odd cycle $\{x_i\}$.

Now consider any two distinct vertices $A$, $B$ of $\mu(G)$.

**Case 1. $A \leftrightarrow B$ is an edge in $\mu(G)$.** It suffices to show that there is a Hamiltonian cycle containing $A \leftrightarrow B$. By symmetry, we may assume $A \leftrightarrow B$ is one of $x_1 \leftrightarrow x_n$, $x_1 \leftrightarrow y_n$, or $y_1 \leftrightarrow z$. In all cases cycle $C$ above contains the edge.

**Case 2. $A \leftrightarrow B$ is not an edge of $\mu(G)$.**

**Case 2.1. $\{A, B\} \subset X$; say $x_i$, $x_j$.** Without loss of generality, assume $i = 1$. The same proof works for both even and odd $j$. Zigzag up from $x_1$ until $y_{j-1}$ is reached (when reaching the end, carry on at the beginning in the obvious way). Then jump via $z$ to $y_n$ and zigzag left (past $y_1$ if necessary) until $x_j$ is reached. Formally: $x_1 \leftrightarrow y_2 \leftrightarrow x_3 \leftrightarrow \cdots \leftrightarrow y_{j-1} \leftrightarrow z \leftrightarrow y_n \leftrightarrow x_{n-1} \leftrightarrow y_{n-2} \leftrightarrow \cdots \leftrightarrow x_j$. Figure 13 shows how this works when $j$ is even, followed by the odd $j$ case.

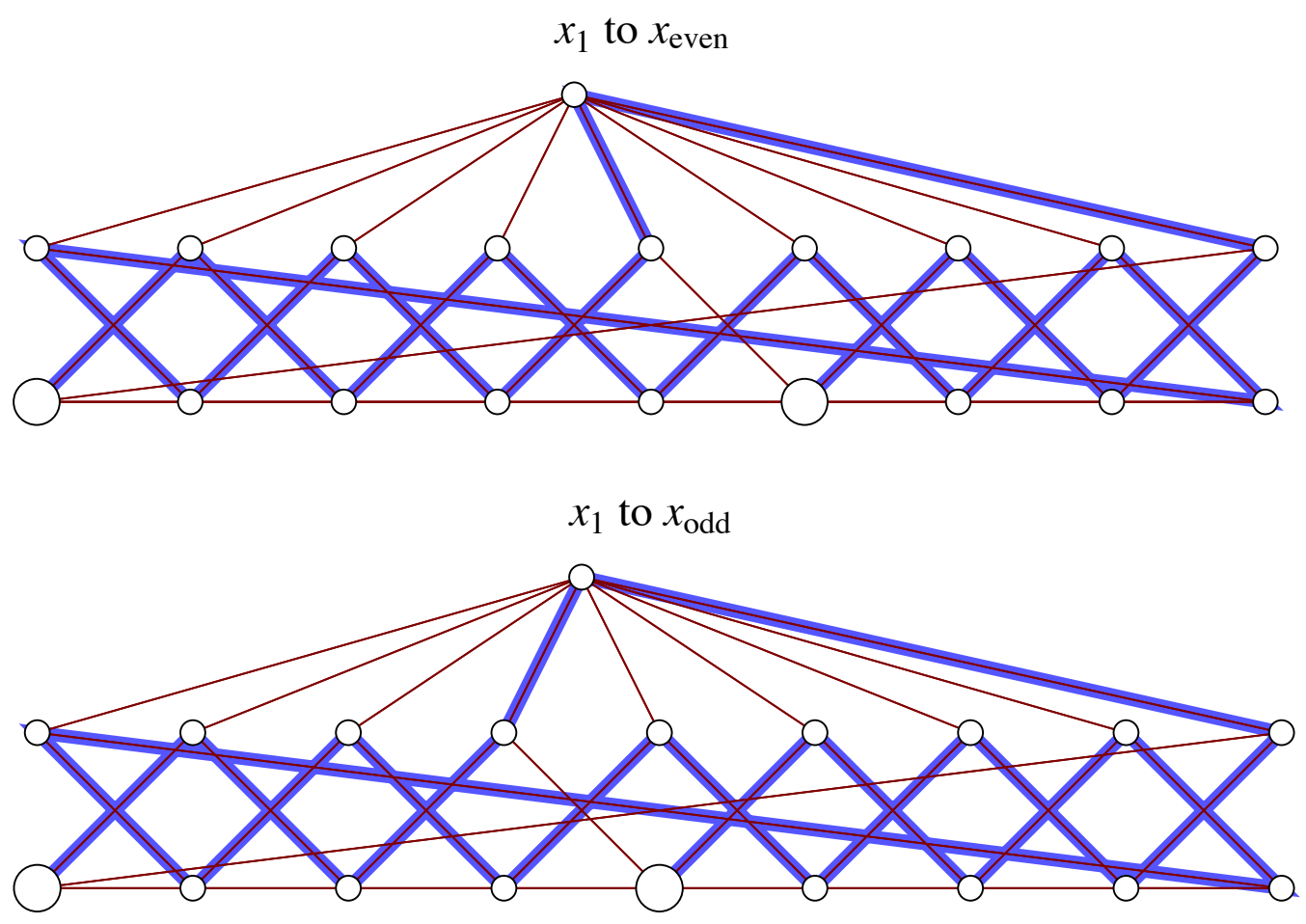

Figure 13. Typical Hamiltonian paths from $X$ to $X$.

**Case 2.2.** ***A*** $\in$ ***X*** **and** ***B*** $\in$ ***Y*****.** Assume $A = x_1$ and $B = y_j$. Then $x_1 \leftrightarrow x_j$ is a nonedge in $G$. If $j$ is even: Zigzag up to $x_{j-1} \leftrightarrow x_j$ then zigzag to $y_n \leftrightarrow z \leftrightarrow y_{j-1}$ and zigzag left and through $y_1$ to the target $x_j$, as in Figure 14.

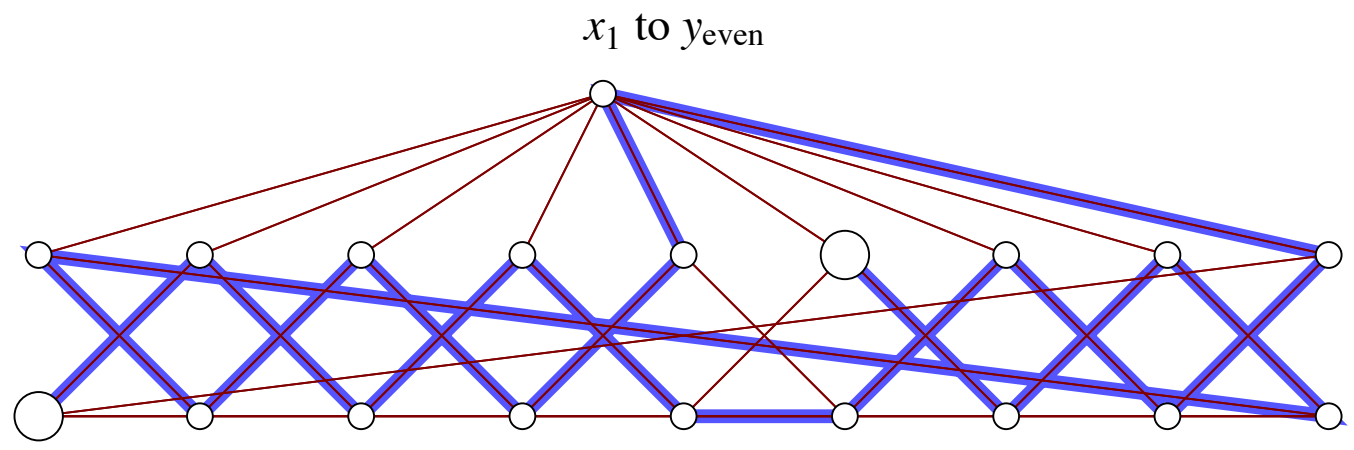

Figure 14. A typical Hamiltonian path from $X$ to an even vertex in $Y$.

If $j$ is odd, zigzag up to $x_j$ then left to $x_{j-1}$ and down through $y_1$ to $y_{j+1}$, then up to $z$ and $y_n$ and zigzag down to the finish at $y_j$, as in Figure 15. This works fine even if $j = n$.

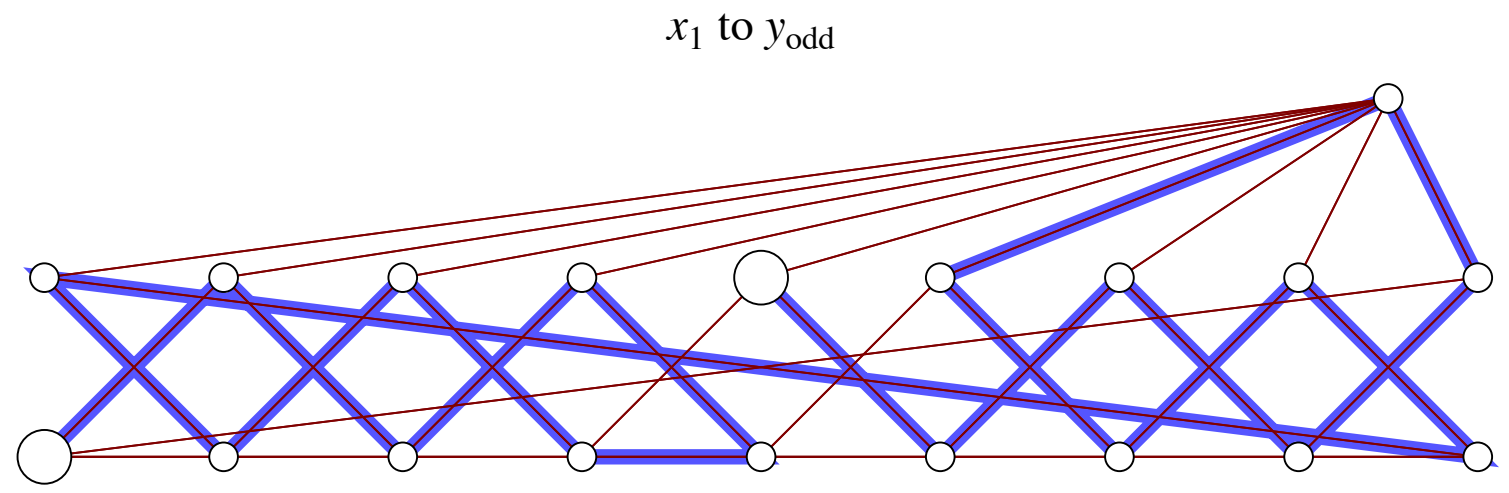

Figure 15. A typical Hamiltonian path from $X$ to an odd vertex in $Y$.

**Case 2.3** ***A*** $\in$ ***X*** **and** ***B*** **=** ***z*****.** Assume $A = x_1$. Zigzag through to $y_n$ and then finish up at $z$:
$x_1 \leftrightarrow y_2 \leftrightarrow x_3 \leftrightarrow y_4 \leftrightarrow x_5 \leftrightarrow y_6 \leftrightarrow \ldots \leftrightarrow x_n \leftrightarrow y_1 \leftrightarrow x_2 \leftrightarrow y_3 \leftrightarrow x_4 \leftrightarrow y_5 \leftrightarrow \ldots \leftrightarrow y_n \leftrightarrow z$; see Figure 16.

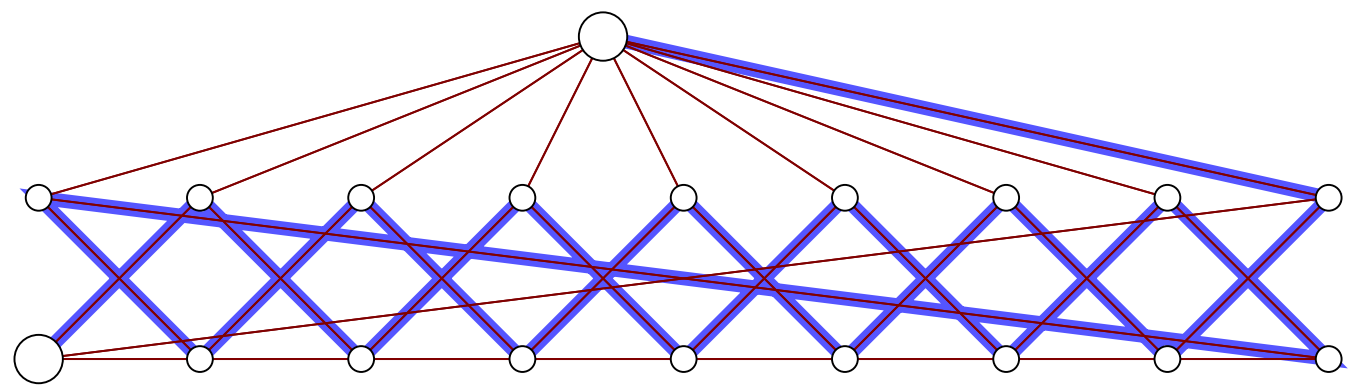

Figure 16. A typical Hamiltonian path from $X$ to $z$.

**Case 2.4** $\{A, B\} \subset Y$. Assume first that $A = y_1$ and $B = y_j$, with $j$ even.

```
{10, 2, 12, 4, 14, 19, 18, 1, 11, 3, 13, 5, 6, 16, 8, 9, 17, 7, 15}
```

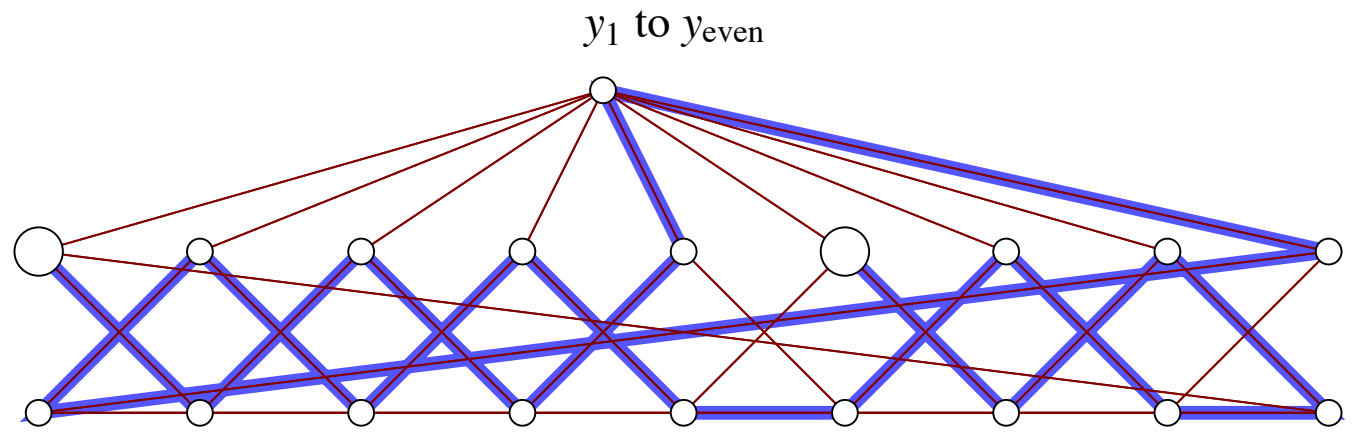

Figure 17. A typical Hamiltonian path from a vertex in $Y$ to a vertex of different parity in $Y$.

Zigzag up from $y_1$ to $y_{j-1}$, then up to $z$ and down to $y_n$, then back to $x_1$ and zigzag up to $x_{j-1}$, then to $x_j$ and zigzag up to $x_{n-1}$, then $x_n$ and zigzag down to $y_j$. See Figure 17. Formally:

$$y_1 \leftrightarrow x_2 \leftrightarrow y_3 \leftrightarrow \cdots \leftrightarrow y_{j-1} \leftrightarrow z \leftrightarrow y_n \leftrightarrow x_1 \leftrightarrow y_2 \leftrightarrow$$
$$x_3 \leftrightarrow \cdots \leftrightarrow x_{j-1} \leftrightarrow x_j \leftrightarrow y_{j+1} \leftrightarrow x_{j+2} \leftrightarrow \cdots \leftrightarrow x_{n-1} \leftrightarrow x_n \leftrightarrow y_{n-1} \leftrightarrow x_{n-2} \cdots \leftrightarrow y_j$$

Finally assume $j$ is odd. Zigzag from $y_1$ to $x_{j-1}$, then back to $x_{j-2}$ and zigzag down to $x_1$ and up to $y_n$, then up to $z$, down to $y_{j-1}$, and zigzag up to $x_n$, then back to $x_{n-1}$ and zigzag to $y_j$; see Figure 18. This completes the proof. □

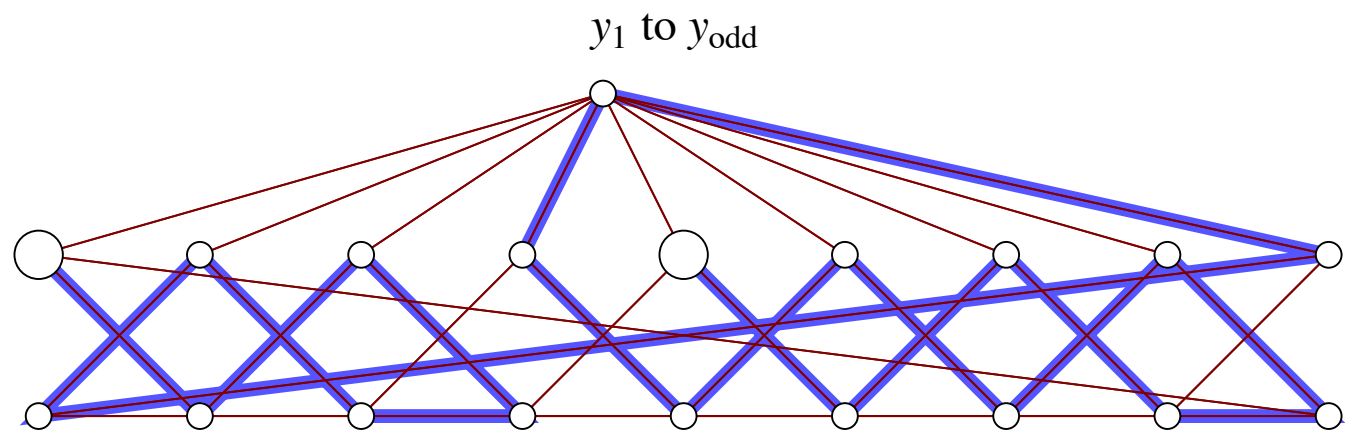

Figure 18. A typical Hamiltonian path from a vertex in $Y$ to a vertex of the same parity in $Y$.

The proof technique does not work directly to settle the even case. But we have barely used the edges in $G$; because $G$ is assumed HC, there might be a way to use more of those edges to extend Corollary 14 to all graphs, proving Conjecture 6. Computation also leads to a conjecture about edge coloring. All the Mycielski graphs (except $M_3$) are class 1 because of Fournier's theorem (Section 1); for $M_4$ and beyond, $z$ is the unique vertex of maximum degree. But perhaps much more is true: computation supports the following conjecture.

**Conjecture 7.** For any graph $G$ other than $K_2$, $\mu(G)$ is class 1.

## 6. Keller Graphs

The *Keller graph* $G_d$ of dimension $d$ is defined as follows [DELSMW, W]: the $4^d$ vertices are all $d$-tuples from $\{0, 1, 2, 3\}$. Two tuples form an edge if they differ in at least two coordinates and if in at least one coordinate the difference of the entries is 2 (mod 4). We ignore $G_1$, which simply consists of four isolated points. These graphs are vertex transitive and therefore regular; it is easy to work out the degree of $G_d$, which is $4^d - 3^d - d$. The graph $G_2$ is also known as the Clebsch graph. The Keller graphs play a critical role in the Keller conjecture [DELSMW], which, in its unrestricted form, states that any tiling of $\mathbb{R}^d$ by unit cubes contains two cubes that meet face-to-face. This conjecture is closely related to $\omega(G_d)$. The value of $\omega(G_d)$ is known for all $d$: when $d \geq 8$, $\omega(G_d) = 2^d$, while $\omega(G_d) < 2^d$ for $d \leq 7$ (Table 2; see [DELSMW]). These $\omega$ values imply that the Keller conjecture with the restriction that all cube centers involve only integers or half-integers is true for $d \leq 7$ and false for $d \geq 8$. The unrestricted Keller conjecture is known to be true for $d \leq 6$ and false for $d \geq 8$, but is unresolved in $\mathbb{R}^7$.

Note that $G_d$ always admits a $2^d$-vertex coloring, defined this way: There are $2^d$ vertices using only 0s and 2s; they each receive a distinct color. Give any other vertex $(v_i)$ the same color as $\left(2\left\lfloor \frac{v_i}{2} \right\rfloor\right)$; the "differ by 2" condition is never satisfied by both vertices $(v_i)$ and $\left(2\left\lfloor \frac{v_i}{2} \right\rfloor\right)$. Therefore $\chi(G_d) \le 2^d$ (proved independently by Fung [F] and Debroni et al [DELSMW]). This coloring is also implicit in the proof of Theorem 17 below: the 0-and-2 set is the diagonal of the array shown. We will show in Corollary 20 that this coloring is optimal for all $d$.

A classic theorem of Dirac [D] states that a graph with minimum degree greater or equal to $n_v/2$ is Hamiltonian; this applies to $G_d$ when $d \ge 3$. We can give an explicit Hamiltonian cycle for all Keller graphs.

**Theorem 17.** All Keller graphs are Hamiltonian.

**Proof.** For $G_2$, a Hamiltonian cycle is (00, 23, 01, 20, 02, 21, 03, 22, 10, 33, 11, 30, 12, 31, 13, 32). This ordering alternates 0 and 2 in the first coordinate for the first half, and then 1 and 3. And in the second coordinate, the leading 0s and 1s are matched, in order, with 0, 1, 2, 3, and the 2s and 3s with 3, 0, 1, 2. For larger $d$, just append vectors to the scheme for $G_2$, thus: $(00X, 23X, 01X, \ldots, 13X, 32X, 00Y, 23Y, \ldots)$, where $X, Y, \ldots$ exhaust all tuples in $\mathbb{Z}_4^{d-2}$. This repetition still yields a cycle and because all vertices are struck, it is Hamiltonian. □

Another classic result [O] states that if the minimum degree of $G$ is greater than or equal to $\frac{1}{2}(n_v + 1)$, then $G$ is Hamilton-connected. The condition holds for $G_d$ when $d \ge 3$ and a simple computation using an algorithm described in [DW] verifies that $G_2$ is Hamilton-connected, so all Keller graphs are HC.

The Keller graphs are vertex-transitive and so provide an infinite family of examples for the conjecture in [DW] that vertex-transitive, Hamiltonian graphs—except cycles and the dodecahedral graph—are HC. Computations also support the conjecture that Keller graphs have *Hamiltonian decompositions* (meaning that the edges can be partitioned into disjoint Hamiltonian cycles, plus a perfect matching if the degree is odd; see Fig. 19 for such a decomposition of $G_2$). We found Hamiltonian decompositions up through $G_6$ and conjecture that they exist for all Keller graphs. Table 1 shows such a decomposition for $G_3$: 17 Hamiltonian cycles.

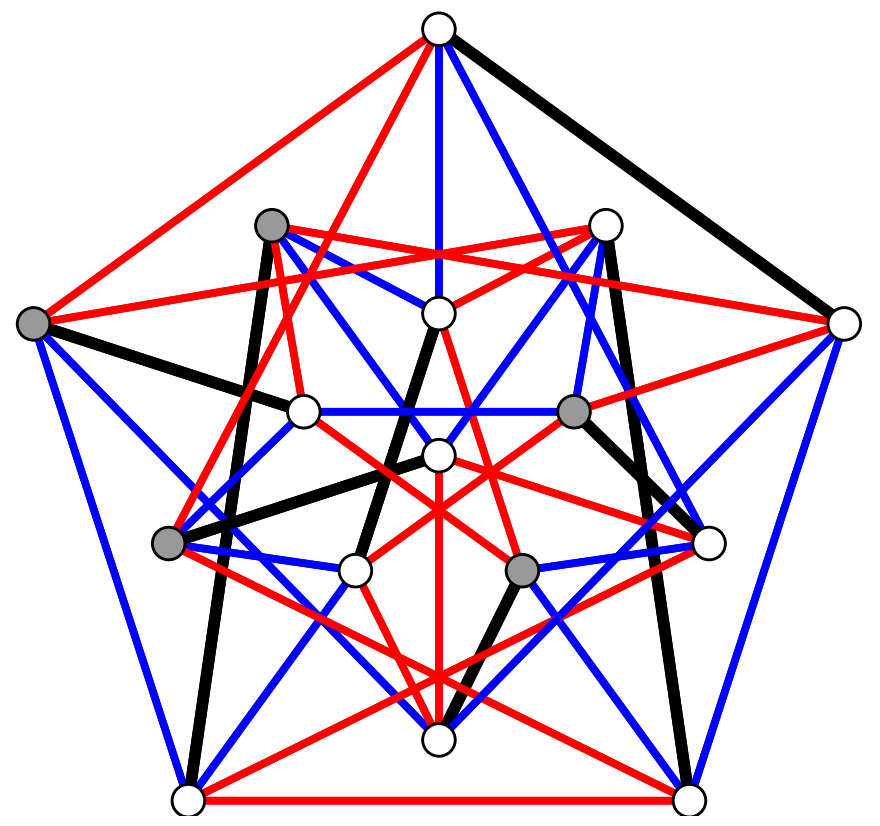

Figure 19. A Hamiltonian decomposition of $G_2$, also known as the Clebsch graph: two Hamiltonian cycles (red, blue) and one perfect matching (black). The gray vertices are a maximum independent set.

| 1 | 2 | 3 | 4 | 5 | 6 | 7 | 8 | 9 | 10 | 11 | 12 | 13 | 14 | 15 | 16 | 17 |
|---|---|---|---|---|---|---|---|---|---|---|---|---|---|---|---|---|
| 0 | 0 | 0 | 0 | 0 | 0 | 0 | 0 | 0 | 0 | 0 | 0 | 0 | 0 | 0 | 0 | 0 |
| 41 | 27 | 54 | 26 | 46 | 30 | 45 | 9 | 22 | 58 | 43 | 25 | 24 | 44 | 59 | 37 | 62 |
| 2 | 9 | 32 | 62 | 60 | 36 | 19 | 31 | 31 | 19 | 53 | 49 | 19 | 1 | 50 | 35 | 22 |
| 45 | 50 | 27 | 55 | 38 | 54 | 8 | 17 | 52 | 48 | 44 | 58 | 33 | 11 | 25 | 1 | 61 |
| 21 | 42 | 60 | 20 | 20 | 8 | 49 | 62 | 46 | 59 | 52 | 12 | 23 | 16 | 63 | 51 | 7 |
| 46 | 48 | 14 | 49 | 53 | 2 | 30 | 29 | 15 | 32 | 62 | 18 | 28 | 40 | 16 | 43 | 21 |
| 16 | 17 | 56 | 15 | 14 | 32 | 21 | 50 | 41 | 11 | 6 | 44 | 51 | 19 | 41 | 50 | 55 |
| 53 | 41 | 42 | 52 | 16 | 56 | 15 | 56 | 5 | 29 | 12 | 14 | 21 | 26 | 3 | 24 | 16 |
| 7 | 34 | 1 | 19 | 51 | 22 | 39 | 27 | 3 | 22 | 40 | 39 | 57 | 63 | 17 | 46 | 58 |
| 28 | 59 | 40 | 25 | 29 | 44 | 33 | 49 | 34 | 40 | 17 | 46 | 43 | 49 | 15 | 32 | 27 |
| 62 | 49 | 33 | 51 | 39 | 9 | 12 | 7 | 40 | 30 | 56 | 22 | 9 | 21 | 45 | 10 | 17 |
| 8 | 23 | 26 | 9 | 12 | 7 | 7 | 38 | 49 | 24 | 29 | 51 | 2 | 52 | 31 | 12 | 35 |
| 39 | 25 | 57 | 39 | 2 | 42 | 40 | 32 | 56 | 2 | 36 | 31 | 59 | 27 | 36 | 21 | 59 |
| 17 | 16 | 20 | 13 | 33 | 3 | 50 | 8 | 30 | 26 | 5 | 21 | 22 | 5 | 9 | 3 | 41 |
| 47 | 62 | 29 | 42 | 25 | 58 | 22 | 18 | 4 | 54 | 42 | 59 | 16 | 39 | 46 | 8 | 19 |
| 1 | 7 | 52 | 33 | 52 | 23 | 4 | 52 | 18 | 27 | 51 | 61 | 6 | 60 | 55 | 47 | 54 |
| 27 | 43 | 12 | 47 | 58 | 29 | 43 | 45 | 42 | 3 | 18 | 35 | 56 | 31 | 28 | 2 | 18 |
| 48 | 13 | 34 | 3 | 17 | 5 | 61 | 54 | 19 | 13 | 24 | 4 | 34 | 4 | 18 | 16 | 53 |
| 25 | 46 | 10 | 61 | 42 | 46 | 52 | 47 | 57 | 51 | 35 | 42 | 13 | 46 | 60 | 59 | 12 |
| 32 | 37 | 45 | 16 | 35 | 6 | 34 | 22 | 16 | 61 | 9 | 16 | 54 | 7 | 22 | 45 | 50 |
| 42 | 14 | 22 | 43 | 41 | 52 | 14 | 13 | 26 | 31 | 38 | 56 | 46 | 57 | 53 | 5 | 28 |
| 24 | 4 | 58 | 45 | 18 | 63 | 42 | 63 | 28 | 53 | 61 | 1 | 3 | 29 | 47 | 34 | 49 |
| 10 | 63 | 18 | 14 | 9 | 17 | 9 | 5 | 61 | 46 | 25 | 10 | 53 | 58 | 6 | 9 | 43 |
| 51 | 53 | 48 | 41 | 59 | 60 | 17 | 35 | 27 | 10 | 3 | 43 | 17 | 32 | 30 | 19 | 37 |
| 40 | 26 | 24 | 31 | 26 | 21 | 54 | 15 | 29 | 4 | 57 | 17 | 10 | 22 | 39 | 56 | 31 |
| 13 | 36 | 63 | 7 | 56 | 39 | 25 | 1 | 59 | 38 | 15 | 24 | 37 | 55 | 11 | 23 | 48 |
| 37 | 18 | 36 | 37 | 28 | 28 | 58 | 25 | 33 | 28 | 20 | 34 | 4 | 61 | 42 | 30 | 26 |
| 19 | 32 | 50 | 59 | 6 | 59 | 36 | 11 | 58 | 34 | 54 | 15 | 60 | 17 | 2 | 58 | 60 |
| 29 | 15 | 16 | 23 | 13 | 20 | 3 | 53 | 2 | 1 | 16 | 33 | 25 | 51 | 40 | 48 | 42 |
| 43 | 5 | 9 | 63 | 31 | 50 | 26 | 59 | 60 | 43 | 30 | 6 | 7 | 33 | 58 | 22 | 6 |
| 3 | 40 | 47 | 6 | 50 | 57 | 44 | 19 | 23 | 18 | 48 | 48 | 32 | 43 | 8 | 52 | 20 |
| 33 | 3 | 5 | 36 | 21 | 24 | 2 | 60 | 13 | 57 | 46 | 57 | 12 | 48 | 44 | 42 | 63 |
| 4 | 29 | 30 | 10 | 27 | 61 | 35 | 20 | 32 | 35 | 23 | 63 | 38 | 8 | 62 | 49 | 33 |
| 34 | 6 | 55 | 60 | 63 | 15 | 56 | 61 | 50 | 49 | 37 | 37 | 5 | 38 | 56 | 26 | 9 |
| 26 | 8 | 19 | 58 | 54 | 25 | 51 | 23 | 11 | 41 | 28 | 47 | 31 | 3 | 20 | 17 | 32 |
| 50 | 22 | 49 | 34 | 15 | 34 | 53 | 57 | 47 | 12 | 4 | 52 | 39 | 59 | 51 | 11 | 5 |
| 52 | 28 | 3 | 27 | 22 | 43 | 27 | 55 | 38 | 62 | 58 | 28 | 49 | 25 | 27 | 38 | 51 |
| 23 | 57 | 44 | 50 | 36 | 12 | 57 | 33 | 24 | 21 | 31 | 36 | 47 | 18 | 35 | 15 | 23 |
| 9 | 51 | 7 | 4 | 11 | 26 | 1 | 24 | 51 | 56 | 59 | 8 | 55 | 62 | 7 | 53 | 45 |
| 49 | 58 | 41 | 32 | 49 | 51 | 23 | 3 | 37 | 25 | 1 | 41 | 18 | 53 | 29 | 25 | 1 |
| 31 | 24 | 23 | 30 | 55 | 45 | 47 | 37 | 8 | 39 | 8 | 27 | 27 | 24 | 49 | 55 | 36 |
| 56 | 60 | 17 | 53 | 45 | 4 | 20 | 44 | 45 | 63 | 26 | 45 | 62 | 6 | 24 | 27 | 2 |
| 18 | 54 | 59 | 29 | 3 | 41 | 46 | 42 | 53 | 7 | 55 | 7 | 40 | 15 | 32 | 13 | 25 |
| 59 | 61 | 13 | 38 | 10 | 55 | 28 | 28 | 39 | 14 | 13 | 13 | 48 | 9 | 43 | 36 | 47 |
| 30 | 47 | 21 | 1 | 48 | 37 | 63 | 58 | 6 | 5 | 19 | 20 | 41 | 37 | 19 | 7 | 29 |
| 54 | 12 | 28 | 46 | 23 | 27 | 18 | 21 | 44 | 33 | 63 | 26 | 1 | 2 | 5 | 60 | 4 |
| 14 | 55 | 53 | 56 | 44 | 16 | 10 | 43 | 10 | 8 | 21 | 32 | 58 | 28 | 23 | 29 | 15 |
| 38 | 1 | 35 | 2 | 5 | 49 | 41 | 2 | 20 | 16 | 47 | 3 | 44 | 10 | 12 | 63 | 40 |
| 63 | 39 | 25 | 57 | 43 | 18 | 13 | 48 | 48 | 52 | 7 | 9 | 50 | 56 | 37 | 41 | 46 |
| 22 | 10 | 31 | 17 | 8 | 11 | 11 | 14 | 9 | 17 | 33 | 40 | 8 | 54 | 1 | 39 | 11 |
| 57 | 19 | 38 | 8 | 34 | 33 | 48 | 36 | 63 | 50 | 27 | 54 | 42 | 23 | 26 | 44 | 34 |
| 11 | 21 | 62 | 35 | 7 | 13 | 29 | 46 | 1 | 23 | 2 | 29 | 36 | 14 | 61 | 4 | 44 |
| 44 | 11 | 39 | 12 | 30 | 53 | 37 | 12 | 7 | 15 | 39 | 55 | 15 | 35 | 21 | 62 | 13 |
| 15 | 35 | 61 | 22 | 37 | 19 | 62 | 30 | 17 | 42 | 45 | 5 | 29 | 13 | 14 | 20 | 52 |
| 55 | 45 | 4 | 24 | 61 | 1 | 32 | 51 | 55 | 20 | 11 | 62 | 35 | 47 | 33 | 18 | 24 |
| 60 | 38 | 2 | 54 | 19 | 31 | 6 | 41 | 14 | 44 | 41 | 23 | 26 | 41 | 10 | 61 | 14 |
| 6 | 44 | 46 | 44 | 62 | 62 | 31 | 6 | 43 | 55 | 32 | 53 | 52 | 50 | 52 | 6 | 8 |
| 35 | 30 | 8 | 21 | 24 | 48 | 55 | 34 | 25 | 9 | 14 | 60 | 14 | 30 | 38 | 40 | 30 |
| 58 | 52 | 51 | 48 | 1 | 38 | 24 | 16 | 62 | 45 | 22 | 30 | 20 | 20 | 13 | 14 | 57 |
| 20 | 2 | 11 | 28 | 32 | 40 | 59 | 10 | 36 | 6 | 49 | 38 | 45 | 34 | 4 | 28 | 39 |
| 12 | 20 | 37 | 5 | 57 | 35 | 5 | 40 | 12 | 37 | 10 | 2 | 63 | 36 | 54 | 54 | 3 |
| 5 | 31 | 15 | 11 | 47 | 10 | 60 | 26 | 54 | 60 | 50 | 11 | 30 | 45 | 48 | 31 | 56 |
| 61 | 33 | 43 | 40 | 4 | 47 | 16 | 4 | 21 | 36 | 60 | 19 | 61 | 12 | 34 | 57 | 38 |
| 36 | 56 | 6 | 18 | 40 | 14 | 38 | 39 | 35 | 47 | 34 | 50 | 11 | 42 | 57 | 33 | 10 |

Table 1. A decomposition of $G_3$ into 17 Hamiltonian cycles. The vertices are encoded, using base 4, by integers from 0 to 63. Because $\Delta = 34$, there are 17 Hamiltonian cycles.

**Conjecture 8.** All Keller graphs have a Hamiltonian decomposition.

Conjecture 8 is related to deep work of Kühn et al [KO, CKLOT]. Theorem 1.7 of [KO] implies that for sufficiently large odd $d$, there is a Hamilton decomposition of $G_d$, while the improvement in [CKLOT, Thm. 1.1.3] handles the even case too. So for sufficiently large $d$, $G_d$ is known to have a Hamiltonian decomposition.

Our algorithm for finding these decompositions starts with the simple idea of trying random class-1 colorings (obtained by using the methods of Vizing and Kempe on a random permutation of the graph) and checking to see if the color-sets, which are matchings, can be paired up to form the desired cycles; the pairing is generally done using the classic blossom algorithm of Edmonds. A more sophisticated approach is needed for large cases such as $G_5$ and $G_6$. We again start with a class-1 coloring, but then apply Kempe switches in the hope of obtaining pairs of matchings that link to form cycles. The heuristic used to decide which Kempe switches to make is a scoring function that compares the number of Hamiltonian cycles obtainable by pairing up matchings (primary key), the minimum number of cycles a pair of remaining matchings produces (secondary key), and the total number of cycles that all pairs of remaining matchings produce (tertiary key). Additionally, a Hamiltonian decomposition of a significant part of the edges of $G_6$ can be effectively constructed from a decomposition of $G_5$. Therefore, to find a decomposition of $G_6$, we first find one for $G_5$. We then apply the local search method using the heuristic function described above, but only to the subgraph of $G_6$ consisting of the uncovered edges.

Although it seemed plausible that connected, vertex-transitive graphs always have Hamiltonian decompositions (excluding a few small examples), that was recently shown by Bryant and Dean [BD] to be false. An even stronger property is that of having a perfect 1-factorization: a collection of matchings such that *any* two form a Hamiltonian cycle. That is a much more difficult subject—it is unresolved even for complete graphs—and all we can say is that an exhaustive search established that $G_2$ does not have a perfect 1-factorization. In the other direction, a weaker conjecture than the false one just mentioned is that all vertex-transitive graphs with even order are class 1 except the Petersen graph and the triangle-replaced Petersen graph; no counterexample is known.

Note that an even-order graph with a Hamiltonian decomposition is necessarily class 1. One can show that all Keller graphs are class 1 by explicit computation up to $G_6$ and then calling on a famous theorem of Chetwynd and Hilton [CW; see also SSTF, Thm. 4.17] for the rest; their theorem applies to graphs for which $\Delta > \frac{1}{2}\left(\sqrt{7} - 1\right) n_v$, which holds for $G_7$ and beyond. But in fact there is a uniform and constructive way to present class-1 colorings of all $G_d$, which we now describe. Note that this result also follows from the class-1 conjecture of the preceding paragraph.

**Theorem 18.** All Keller graphs are class 1.

**Proof.** All arithmetic here is mod 4. Call a vertex—a $d$-tuple—*even* if all entries are even; otherwise odd. The class-1 coloring can be constructed explicitly as follows. Define the color set $S$ to consist of all vertices whose coordinates have at least one 2, but excluding the $d$ vectors consisting of just $d - 1$ 0s and one 2. This set is a type of kernel: the set of all differences $u - v$ for edges $v \longleftrightarrow u$. This set satisfies (1) $|S| = \Delta$, and (2) for each vertex $v$, its neighbors are $v + S$. Partition $S$ into its even vectors, $S_0$, and its odd ones, $S_1$.

For each $s \in S_1$, define an equivalence relation $\sim_s$ on the vertices: $u \sim_s v$ iff $u - v$ is a multiple of $s$. Each equivalence class has the form $\{v, v + s, v + 2s, v + 3s\}$; because $s$ is odd, each such class has four distinct elements. Note that the collection of classes for $s$ is identical to the collection of classes for $-s$. For each $s \in S_1$, define a choice set $C_s$ for the equivalence classes; use the lexicographically first vector in each class. Then $C_s = C_{-s}$.

Define the edge coloring as follows (see Fig. 20). For each even color $s \in S_0$, use it for all edges $v \longleftrightarrow v + s$. Each $s$ colors $\frac{1}{2} n_v$ edges because $v \longleftrightarrow v \pm s$ both get color $s$, but are the same edge (because $s = -s$). So in all, this colors $\frac{1}{2} n_v \, |S_0|$ edges. For each odd color $s \in S_1$, use it for the edges $v \longleftrightarrow v + s$ and $v + 2s \longleftrightarrow v + 3s$, but, in both cases, only for vertices $v \in C_s$. Because $|C_s| = \frac{1}{4} n_v$, each color applies to $\frac{1}{2} n_v$ edges, and so the odd colors taken together color $\frac{1}{2} n_v \, |S_1|$ edges. Thus the number of edges that are colored is $\frac{1}{2} n_v \, (|S_0| + |S_1|) = \frac{1}{2} n_v \, |S| = \frac{1}{2} n_v \, \Delta = n_e$, the total number of edges.

*Claim.* Every edge receives only one color. *Proof.* Given edge $u \longleftrightarrow w$, let $s = w - u$. If $s$ is even, then $s = -s$ and this easily yields the claim. The odd case is more delicate. Suppose $u \longleftrightarrow w$ is assigned color $s$; then $s = \pm(w - u)$. We may assume $s = w - u$. If the edge is assigned another color distinct from $s$, that color must therefore be $-s$. Now $u$ and $w$ are equivalent under both relations $\sim_s$ and $\sim_{-s}$. And the class representatives from $C_s$ and $C_{-s}$ agree. This means that $u \longleftrightarrow w$ must be one of the edges $\{v \longleftrightarrow v + s,\; v + 2s \longleftrightarrow v + 3s\}$ and also one of the edges $\{v \longleftrightarrow v - s,\; v - 2s \longleftrightarrow v - 3s\}$. But the latter set equals $\{v + 3s \longleftrightarrow v,\; v + s \longleftrightarrow v + 2s\}$, which is disjoint from the first pair.

The claim and the fact that $n_e$ edges are colored means that every edge receives a color. So it remains only to show that the coloring is proper. Suppose not. Then we have edges $u \longleftrightarrow w$ and $u \longleftrightarrow y$ receiving the same color $s$. If $s$ is even, this is not possible because the edges would have to be of the form $v \longleftrightarrow v + s$ and $v \longleftrightarrow v - s$, which are equal because $s = -s$. Suppose $s$ is odd and color $s$ is assigned to edge $u \longleftrightarrow w$. If $u \in C_s$, then $w = u + s$; if $u = v + s$ where $v \in C_s$, then $w = v$; if $u = v + 2s$ where $v \in C_s$, then $w = v + 3s$, and if $u = v + 3s$ where $v \in C_s$, then $w = v + 2s$. In all cases there is only one choice for $w$. □

| colors $S$ | size-4 equivalence classes of vertices with edges colored by the odd element of $S$ | | | | | | | |
|---|---|---|---|---|---|---|---|---|
| ↑ 23 | 00–23 | 02–21 | 01–20 | 03–22 | 10–33 | 12–31 | 11–30 | 13–32 |
| 21 | 00–21 | 02–23 | 01–22 | 03–20 | 10–31 | 12–33 | 11–32 | 13–30 |
| 32 | 00–32 | 20–12 | 01–33 | 21–13 | 02–30 | 22–10 | 03–31 | 23–11 |
| odd 12 | 00–12 | 20–32 | 01–13 | 21–33 | 02–10 | 22–30 | 03–11 | 23–31 |
| even ↓ 22 | 00–22 | 01–23 | 02–20 | 03–21 | 10–32 | 11–33 | 12–30 | 13–31 |

Figure 20. The class-1 Keller coloring for the edges of $G_2$, using five colors. The even case has only one entry; $S_0 = \{22\}$. The odd case has four colors and the four equivalence classes of the full vertex set are shown, with the matchings within each class. Note that the classes for $\pm s$ are the same sets (e.g., $s = 12$ and 32).

We can also investigate some familiar parameters for Keller graphs. The standard parameters $\alpha$, $\theta$, $\chi$, and $\chi_{\text{frac}}$ are defined in Section 1. Let $\theta_{\text{frac}}$ be the fractional clique covering number (same as $\chi_{\text{frac}}$ of the complementary graph). Table 2 shows the known results, including results proved here. It is clear that $\alpha(G_d) \geq 2^d$ since the tuples using only 0s and 1s are independent. A larger independent set can exist, but only in $G_2$, as Theorem 19 shows.

| $d$ | all $d$ | 2 | 3 | 4 | 5 | 6 | 7 | $d \geq 8$ |
|---|---|---|---|---|---|---|---|---|
| independence number, $\alpha$ | | 5 | 8 | 16 | 32 | 64 | 128 | $2^d$ |
| chromatic number, $\chi$ | $2^d$ | 4 | 8 | 16 | 32 | 64 | 128 | $2^d$ |
| fractional chromatic number, $\chi_{\text{frac}}$ | | $\frac{16}{5}$ | 8 | 16 | 32 | 64 | 128 | $2^d$ |
| class 1 for edge coloring; $\chi' = \Delta$ | Yes | | | | | | | |
| maximum clique size, $\omega$ | | 2 | 5 | 12 | 28 | 60 | 124 | $2^d$ |
| clique covering number, $\theta$ | | 8 | 13 | 22 | $37 \leq \theta \leq 40$ | $69 \leq \theta \leq 80$ | $133 \leq \theta \leq 160$ | $2^d$ |
| fractional clique covering number, $\theta_{\text{frac}}$ | | 8 | $\frac{64}{5}$ | $\frac{64}{3}$ | $\frac{256}{7}$ | $\frac{1024}{15}$ | $\frac{4096}{31}$ | $2^d$ |
| Hamiltonian, Hamilton-connected | Yes | | | | | | | |
| Hamiltonian decomposition | Conjectured yes | Yes | Yes | Yes | Yes | Yes | ? | ? |
| perfect 1-factorization | | No | ? | ? | ? | ? | ? | ? |
| degree, $\Delta$ | $4^d - 3^d - d$ | 5 | 34 | 171 | 776 | 3361 | 14190 | |

Table 2. Properties of the Keller graphs $G_d$. The number of vertices of $G_d$ is $4^d$ and the edge count is $\frac{1}{2} 4^d \left(4^d - 3^d - d\right)$.

**Theorem 19.** The independence number of $G_d$ is $2^d$, except that $\alpha(G_2) = 5$.

**Proof**. For $d \leq 5$, this was known; direct computational methods work. The anomalous case has maximum independent set $\{(0, 3), (1, 0), (1, 2), (1, 3), (2, 3)\}$; see Figure 19. For $d = 6$ or 7, one can again use computation, but some efficiencies are needed since the graphs are large. The graphs are vertex-transitive, so we may assume the first vertex is in the largest independent set. Thus, if $A$ consists of the first vertex together with its neighbors, we can look at $H_d$, the subgraph of $G_d$ generated by the vertices not in $A$. This is substantially smaller, and we need only show that $\alpha(H_d) = 2^d - 1$. That can be done by standard algorithms for finding independent sets; in *Mathematica* it takes a fraction of a second to show that this is the case for $H_6$ and only a few seconds to do the same for $H_7$.

Now suppose $d \geq 8$. Recall that it is known that $\omega(G_d) = 2^d$ in this case (Mackey [M] for $d = 8$; see [DELSMW, Thm. 4.2] for larger $d$). As in [DELSMW], place the vertex labels in a $2^d \times 2^d$ grid, called the *independence square*. The row position of a tuple is computed by converting 0 or 1 to 0 and also converting 2 or 3 to 1 and then treating the result as a binary number. The column position of a tuple is computed by converting 0 or 3 to 0 and also converting 1 or 2 to 1 and then interpreting this in binary. The array for $G_4$ is shown in Table 3.

| | | | | | | | | | | | | | | | |
|---|---|---|---|---|---|---|---|---|---|---|---|---|---|---|---|
| 0000 | 0001 | 0010 | 0011 | 0100 | 0101 | 0110 | 0111 | 1000 | 1001 | 1010 | 1011 | 1100 | 1101 | 1110 | 1111 |
| 0003 | 0002 | 0013 | 0012 | 0103 | 0102 | 0113 | 0112 | 1003 | 1002 | 1013 | 1012 | 1103 | 1102 | 1113 | 1112 |
| 0030 | 0031 | 0020 | 0021 | 0130 | 0131 | 0120 | 0121 | 1030 | 1031 | 1020 | 1021 | 1130 | 1131 | 1120 | 1121 |
| 0033 | 0032 | 0023 | 0022 | 0133 | 0132 | 0123 | 0122 | 1033 | 1032 | 1023 | 1022 | 1133 | 1132 | 1123 | 1122 |
| 0300 | 0301 | 0310 | 0311 | 0200 | 0201 | 0210 | 0211 | 1300 | 1301 | 1310 | 1311 | 1200 | 1201 | 1210 | 1211 |
| 0303 | 0302 | 0313 | 0312 | 0203 | 0202 | 0213 | 0212 | 1303 | 1302 | 1313 | 1312 | 1203 | 1202 | 1213 | 1212 |
| 0330 | 0331 | 0320 | 0321 | 0230 | 0231 | 0220 | 0221 | 1330 | 1331 | 1320 | 1321 | 1230 | 1231 | 1220 | 1221 |
| 0333 | 0332 | 0323 | 0322 | 0233 | 0232 | 0223 | 0222 | 1333 | 1332 | 1323 | 1322 | 1233 | 1232 | 1223 | 1222 |
| 3000 | 3001 | 3010 | 3011 | 3100 | 3101 | 3110 | 3111 | 2000 | 2001 | 2010 | 2011 | 2100 | 2101 | 2110 | 2111 |
| 3003 | 3002 | 3013 | 3012 | 3103 | 3102 | 3113 | 3112 | 2003 | 2002 | 2013 | 2012 | 2103 | 2102 | 2113 | 2112 |
| 3030 | 3031 | 3020 | 3021 | 3130 | 3131 | 3120 | 3121 | 2030 | 2031 | 2020 | 2021 | 2130 | 2131 | 2120 | 2121 |
| 3033 | 3032 | 3023 | 3022 | 3133 | 3132 | 3123 | 3122 | 2033 | 2032 | 2023 | 2022 | 2133 | 2132 | 2123 | 2122 |
| 3300 | 3301 | 3310 | 3311 | 3200 | 3201 | 3210 | 3211 | 2300 | 2301 | 2310 | 2311 | 2200 | 2201 | 2210 | 2211 |
| 3303 | 3302 | 3313 | 3312 | 3203 | 3202 | 3213 | 3212 | 2303 | 2302 | 2313 | 2312 | 2203 | 2202 | 2213 | 2212 |
| 3330 | 3331 | 3320 | 3321 | 3230 | 3231 | 3220 | 3221 | 2330 | 2331 | 2320 | 2321 | 2230 | 2231 | 2220 | 2221 |
| 3333 | 3332 | 3323 | 3322 | 3233 | 3232 | 3223 | 3222 | 2333 | 2332 | 2323 | 2322 | 2233 | 2232 | 2223 | 2222 |

Table 3. The independence square for $G_4$: each row and each column is an independent set.

The tuples in the same row of the square form an independent set because, in each digit, the value is always either 0 or 1, or it is 2 or 3. Therefore there is no position where the difference is 2 (mod 4). Similarly, the tuples in a column form an independent set. The independence square also proves that $\chi(G_d) \leq 2^d$ for any $d$.

Let $X$ be a clique of order $2^d$ in $G_d$. It has exactly one entry per row in the independence square. Given a $d$-digit bit-string $b$, use it to define an associated automorphism of the graph. In the positions where $b$ has 0, leave the corresponding position of all the vertex entries alone. In the places where $b$ has 1, do the following in those positions: switch 1 and 0, and switch 2 and 3. This preserves adjacency because positions that are different in value are still different and positions that differed by 2 (mod 4) still differ by 2 (mod 4).

For example: if the bit string was 0011, then the first two columns stay the same and the last two have the swaps: 0213 becomes 0202. The complete action of the automorphism on $G_3$ using the bit-string 001 is shown in Table 4.

| | | | | | | | |
|---|---|---|---|---|---|---|---|
| 000 | 001 | 010 | 011 | 100 | 101 | 110 | 111 |
| 003 | 002 | 013 | 012 | 103 | 102 | 113 | 112 |
| 030 | 031 | 020 | 021 | 130 | 131 | 120 | 121 |
| 033 | 032 | 023 | 022 | 133 | 132 | 123 | 122 |
| 300 | 301 | 310 | 311 | 200 | 201 | 210 | 211 |
| 303 | 302 | 313 | 312 | 203 | 202 | 213 | 212 |
| 330 | 331 | 320 | 321 | 230 | 231 | 220 | 221 |
| 333 | 332 | 323 | 322 | 233 | 232 | 223 | 222 |

(a)

| | | | | | | | |
|---|---|---|---|---|---|---|---|
| 001 | 000 | 011 | 010 | 101 | 100 | 111 | 110 |
| 002 | 003 | 012 | 013 | 102 | 103 | 112 | 113 |
| 031 | 030 | 021 | 020 | 131 | 130 | 121 | 120 |
| 032 | 033 | 022 | 023 | 132 | 133 | 122 | 123 |
| 301 | 300 | 311 | 310 | 201 | 200 | 211 | 210 |
| 302 | 303 | 312 | 313 | 202 | 203 | 212 | 213 |
| 331 | 330 | 321 | 320 | 231 | 230 | 221 | 220 |
| 332 | 333 | 322 | 323 | 232 | 233 | 222 | 223 |

(b)

Table 4. (a) The independence square for $G_3$. (b) After application of the automorphism defined by 001.

Note that this automorphism maps each row of the square to itself. The collection of automorphisms that correspond to all 0-1 bit strings will map a $2^d$-clique of the Keller graph to a partitioning of the vertices of the Keller graph into $2^d$ disjoint cliques each of size $2^d$. So $\theta(G_d) = 2^d$ for such graphs.

Since any coloring of the graph can have at most one vertex per clique, for Keller graphs that have a $2^d$ clique (i.e., for $d \geq 8$, which we have assumed), it is not possible to find an independent set of size bigger than $2^d$. □

**Corollary 20.** For all Keller graphs, $\chi(G_d) = 2^d$.

**Proof.** The constructive coloring at the beginning of the section gives $2^d$ as an upper bound. Theorem 19 gives $2^d$ as a lower bound, because $4^d/\alpha(G_d) = 2^d$. □

Fung [F, Cor. 6.7] observed that $\chi(G_d) = 2^d$ for $d \geq 4$, $\chi(G_3) \geq 7$, and $\chi(G_2) \geq 3$. Corollary 20 establishes the validity of $\chi = 2^d$ for all Keller graphs.

The graph $G_2$ is an anomaly, with independence number 5. Because $\chi_{\text{frac}}(G) = \frac{n_v}{\alpha(G)}$ for vertex-transitive graphs (see [SU]), we get the following result.

**Corollary 21.** If $d \geq 3$, then $\chi_{\text{frac}}(G_d) = 2^d$; $\chi_{\text{frac}}(G_2) = 16/5$.

So for $d \geq 8$, we have that each parameter $\alpha$, $\chi$, $\chi_{\text{frac}}$ and $\omega$ equals $2^d$.

Computing $\theta(G_d)$ when $d \leq 7$ is difficult. A general lower bound is $\left\lceil \frac{n_V(H)}{\omega(H)} \right\rceil \leq \theta(H)$, which yields 8, 13, 22, 37, 69, 133, the lower bounds of Table 2 (note also that $\alpha \leq \theta$); the first three are sharp. But for $5 \leq d \leq 7$, we do not know $\theta(G_d)$. For $G_5$, we have only that $37 \leq \theta \leq 40$. Because $G_2$ has no triangles, it is clear that $\theta(G_2) = 8$. A 13-coloring of the complement of $G_3$ is shown in Table 5. Table 6 shows a covering of $G_4$ by 22 cliques, the method for which we will explain shortly. That same method found $\theta(G_5) \leq 40$ (see Table 7). The values of $\theta_{\text{frac}}$ in Table 2 arise from the vertex-transitive formula $\chi_{\text{frac}} = n_v/\alpha$ on the complement, which becomes $n_v/\omega$.

| | | | | | |
|---|---|---|---|---|---|
| 1 | 113 | 130 | 232 | 300 | 312 |
| 2 | 110 | 131 | 212 | 320 | 332 |
| 3 | 100 | 121 | 202 | 310 | 322 |
| 4 | 021 | 033 | 203 | 220 | 301 |
| 5 | 031 | 112 | 133 | 303 | 311 |
| 6 | 012 | 020 | 132 | 213 | 230 |
| 7 | 013 | 032 | 111 | 223 | 231 |
| 8 | 011 | 030 | 201 | 233 | 313 |
| 9 | 010 | 022 | 200 | 221 | 302 |
| 10 | 003 | 101 | 122 | 323 | 331 |
| 11 | 001 | 103 | 120 | 321 | 333 |
| 12 | 000 | 023 | 102 | 210 | 222 |
| 13 | 002 | 123 | 211 | 330 | |

Table 5. A covering of the 64 vertices of $G_3$ by 13 disjoint complete subgraphs.

| | | | | | | | | | | | | |
|---|---|---|---|---|---|---|---|---|---|---|---|---|
| 1 | 0233 | 1003 | 1020 | 1213 | 1221 | 1301 | 2101 | 3033 | 3113 | 3121 | 3312 | 3331 |
| 2 | 0013 | 0030 | 0200 | 1212 | 1220 | 1332 | 2012 | 2020 | 2100 | 3032 | 3202 | 3221 |
| 3 | 0100 | 0132 | 0212 | 0310 | 0333 | 1331 | 2010 | 2033 | 2131 | 2211 | 2223 | 3012 |
| 4 | 0223 | 0303 | 0331 | 1102 | 1121 | 1311 | 2123 | 2313 | 2330 | 3011 | 3103 | 3131 |
| 5 | 0101 | 1022 | 1103 | 1120 | 1302 | 1330 | 2222 | 3013 | 3021 | 3203 | 3220 | 3301 |
| 6 | 0001 | 0033 | 0113 | 1021 | 1211 | 1232 | 2000 | 2023 | 2213 | 3201 | 3233 | 3321 |
| 7 | 0111 | 0123 | 0203 | 0301 | 0322 | 1320 | 2001 | 2022 | 2120 | 2200 | 2232 | 3003 |
| 8 | 0302 | 1111 | 1132 | 1230 | 1310 | 1322 | 2030 | 3010 | 3022 | 3102 | 3200 | 3223 |
| 9 | 0003 | 0020 | 0122 | 0202 | 0230 | 1001 | 2113 | 2121 | 2201 | 2303 | 2320 | 3322 |
| 10 | 0011 | 0103 | 0131 | 0313 | 0330 | 1123 | 2002 | 2021 | 2203 | 2231 | 2323 | 3211 |
| 11 | 0110 | 1033 | 1112 | 1131 | 1313 | 1321 | 2233 | 3002 | 3030 | 3212 | 3231 | 3310 |
| 12 | 0012 | 0031 | 0133 | 0213 | 0221 | 1010 | 2102 | 2130 | 2210 | 2312 | 2331 | 3333 |
| 13 | 0220 | 1011 | 1023 | 1201 | 1222 | 1303 | 2103 | 3020 | 3101 | 3122 | 3300 | 3332 |
| 14 | 0000 | 0032 | 0120 | 0201 | 0222 | 1012 | 2110 | 2133 | 2212 | 2300 | 2332 | 3320 |
| 15 | 0231 | 1000 | 1032 | 1210 | 1233 | 1312 | 2112 | 3031 | 3110 | 3133 | 3311 | 3323 |
| 16 | 0121 | 0311 | 0332 | 1013 | 1101 | 1133 | 2221 | 2301 | 2333 | 3100 | 3123 | 3313 |
| 17 | 0002 | 0021 | 0211 | 1203 | 1231 | 1323 | 2003 | 2031 | 2111 | 3023 | 3213 | 3230 |
| 18 | 0232 | 0312 | 0320 | 1113 | 1130 | 1300 | 2132 | 2302 | 2321 | 3000 | 3112 | 3120 |
| 19 | 0010 | 0022 | 0102 | 1030 | 1200 | 1223 | 2011 | 2032 | 2202 | 3210 | 3222 | 3330 |
| 20 | 0130 | 0300 | 0323 | 1002 | 1110 | 1122 | 2230 | 2310 | 2322 | 3111 | 3132 | 3302 |
| 21 | 0112 | 0321 | 1100 | 1333 | 2122 | 2311 | 3130 | 3303 | | | | |
| 22 | 0023 | 0210 | 1031 | 1202 | 2013 | 2220 | 3001 | 3232 | | | | |

Table 6. A covering of the 256 vertices of $G_4$ by 22 disjoint complete subgraphs.

| | | | | | | | | | | | | | | | | | | | | | | | | | | | | |
|---|---|---|---|---|---|---|---|---|---|---|---|---|---|---|---|---|---|---|---|---|---|---|---|---|---|---|---|---|
| 1 | 0 | 18 | 328 | 346 | 44 | 382 | 401 | 147 | 457 | 203 | 176 | 242 | 424 | 490 | 529 | 595 | 777 | 843 | 624 | 882 | 552 | 810 | 925 | 719 | 753 | 739 | 953 | 939 |
| 2 | 65 | 343 | 25 | 31 | 33 | 51 | 56 | 314 | 393 | 399 | 421 | 439 | 172 | 430 | 533 | 791 | 524 | 542 | 564 | 562 | 620 | 890 | 641 | 899 | 904 | 922 | 932 | 930 |
| 3 | 276 | 22 | 332 | 78 | 117 | 55 | 365 | 303 | 149 | 135 | 477 | 463 | 497 | 163 | 784 | 578 | 628 | 614 | 828 | 814 | 660 | 726 | 908 | 974 | 1013 | 759 | 941 | 687 |
| 4 | 260 | 6 | 348 | 94 | 357 | 295 | 125 | 63 | 389 | 407 | 205 | 223 | 225 | 435 | 848 | 514 | 884 | 870 | 572 | 558 | 964 | 902 | 732 | 670 | 741 | 999 | 701 | 959 |
| 5 | 5 | 263 | 93 | 351 | 356 | 294 | 124 | 62 | 388 | 406 | 204 | 222 | 224 | 434 | 513 | 851 | 869 | 887 | 557 | 575 | 901 | 967 | 669 | 735 | 740 | 998 | 700 | 958 |
| 6 | 257 | 275 | 73 | 91 | 109 | 319 | 128 | 386 | 216 | 474 | 481 | 419 | 249 | 187 | 768 | 834 | 536 | 602 | 609 | 867 | 569 | 827 | 652 | 990 | 992 | 1010 | 680 | 698 |
| 7 | 1 | 19 | 329 | 347 | 45 | 383 | 400 | 146 | 456 | 202 | 177 | 243 | 425 | 491 | 528 | 594 | 776 | 842 | 625 | 883 | 553 | 811 | 924 | 718 | 752 | 738 | 952 | 938 |
| 8 | 4 | 262 | 92 | 350 | 293 | 359 | 61 | 127 | 405 | 391 | 221 | 207 | 433 | 227 | 512 | 850 | 868 | 886 | 556 | 574 | 900 | 966 | 668 | 734 | 997 | 743 | 957 | 703 |
| 9 | 69 | 339 | 269 | 267 | 289 | 307 | 40 | 298 | 157 | 155 | 165 | 183 | 188 | 446 | 517 | 775 | 780 | 798 | 800 | 806 | 616 | 894 | 657 | 915 | 648 | 666 | 688 | 694 |
| 10 | 320 | 326 | 97 | 355 | 376 | 362 | 208 | 214 | 408 | 142 | 245 | 503 | 252 | 238 | 837 | 855 | 844 | 590 | 873 | 879 | 705 | 723 | 984 | 730 | 673 | 951 | 761 | 767 |
| 11 | 80 | 86 | 113 | 371 | 104 | 122 | 448 | 454 | 136 | 414 | 229 | 487 | 492 | 510 | 581 | 599 | 604 | 862 | 637 | 635 | 961 | 979 | 712 | 970 | 677 | 947 | 1005 | 1003 |
| 12 | 341 | 67 | 29 | 27 | 49 | 35 | 312 | 58 | 397 | 395 | 437 | 423 | 428 | 174 | 789 | 535 | 540 | 526 | 560 | 566 | 888 | 622 | 897 | 643 | 920 | 906 | 928 | 934 |
| 13 | 325 | 323 | 352 | 98 | 361 | 379 | 213 | 211 | 141 | 411 | 500 | 246 | 237 | 255 | 852 | 838 | 589 | 847 | 876 | 874 | 720 | 706 | 729 | 987 | 948 | 674 | 764 | 762 |
| 14 | 277 | 23 | 333 | 79 | 52 | 118 | 300 | 366 | 132 | 150 | 460 | 478 | 160 | 498 | 785 | 579 | 629 | 615 | 829 | 815 | 661 | 727 | 909 | 975 | 756 | 1014 | 684 | 942 |
| 15 | 16 | 2 | 344 | 330 | 380 | 46 | 145 | 403 | 201 | 459 | 240 | 178 | 488 | 426 | 593 | 531 | 841 | 779 | 880 | 626 | 808 | 554 | 717 | 927 | 737 | 755 | 937 | 955 |
| 16 | 340 | 66 | 28 | 26 | 48 | 34 | 313 | 59 | 396 | 394 | 436 | 422 | 429 | 175 | 788 | 534 | 541 | 527 | 561 | 567 | 889 | 623 | 896 | 642 | 921 | 907 | 929 | 935 |
| 17 | 81 | 87 | 112 | 370 | 105 | 123 | 449 | 455 | 137 | 415 | 228 | 486 | 493 | 511 | 580 | 598 | 605 | 863 | 636 | 634 | 960 | 978 | 713 | 971 | 676 | 946 | 1004 | 1002 |
| 18 | 273 | 259 | 89 | 75 | 317 | 111 | 384 | 130 | 472 | 218 | 417 | 483 | 185 | 251 | 832 | 770 | 600 | 538 | 865 | 611 | 825 | 571 | 988 | 654 | 1008 | 994 | 696 | 682 |
| 19 | 84 | 82 | 369 | 115 | 120 | 106 | 452 | 450 | 412 | 138 | 485 | 231 | 508 | 494 | 597 | 583 | 860 | 606 | 633 | 639 | 977 | 963 | 968 | 714 | 945 | 679 | 1001 | 1007 |
| 20 | 272 | 258 | 88 | 74 | 316 | 110 | 129 | 387 | 217 | 475 | 416 | 482 | 184 | 250 | 769 | 835 | 537 | 603 | 864 | 610 | 824 | 570 | 653 | 991 | 993 | 1011 | 681 | 699 |
| 21 | 336 | 70 | 264 | 270 | 304 | 290 | 297 | 43 | 152 | 158 | 180 | 166 | 445 | 191 | 772 | 518 | 797 | 783 | 805 | 803 | 893 | 619 | 912 | 658 | 665 | 651 | 693 | 691 |
| 22 | 85 | 83 | 368 | 114 | 121 | 107 | 453 | 451 | 413 | 139 | 484 | 230 | 509 | 495 | 596 | 582 | 861 | 607 | 632 | 638 | 976 | 962 | 969 | 715 | 944 | 678 | 1000 | 1006 |
| 23 | 17 | 3 | 345 | 331 | 381 | 47 | 144 | 402 | 200 | 458 | 241 | 179 | 489 | 427 | 592 | 530 | 840 | 778 | 881 | 627 | 809 | 555 | 716 | 926 | 736 | 754 | 936 | 954 |
| 24 | 337 | 71 | 265 | 271 | 305 | 291 | 296 | 42 | 153 | 159 | 181 | 167 | 444 | 190 | 773 | 519 | 796 | 782 | 804 | 802 | 892 | 618 | 913 | 659 | 664 | 650 | 692 | 690 |
| 25 | 256 | 274 | 72 | 90 | 108 | 318 | 385 | 131 | 473 | 219 | 480 | 418 | 248 | 186 | 833 | 771 | 601 | 539 | 608 | 866 | 568 | 826 | 989 | 655 | 1009 | 995 | 697 | 683 |
| 26 | 21 | 279 | 77 | 335 | 116 | 54 | 364 | 302 | 148 | 134 | 476 | 462 | 496 | 162 | 577 | 787 | 613 | 631 | 813 | 831 | 725 | 663 | 973 | 911 | 1012 | 758 | 940 | 686 |
| 27 | 324 | 322 | 353 | 99 | 360 | 378 | 212 | 210 | 140 | 410 | 501 | 247 | 236 | 254 | 853 | 839 | 588 | 846 | 877 | 875 | 721 | 707 | 728 | 986 | 949 | 675 | 765 | 763 |
| 28 | 261 | 7 | 349 | 95 | 292 | 358 | 60 | 126 | 404 | 390 | 220 | 206 | 432 | 226 | 849 | 515 | 885 | 871 | 573 | 559 | 965 | 903 | 733 | 671 | 996 | 742 | 956 | 702 |
| 29 | 64 | 342 | 24 | 30 | 32 | 50 | 57 | 315 | 392 | 398 | 420 | 438 | 173 | 431 | 532 | 790 | 525 | 543 | 565 | 563 | 621 | 891 | 640 | 898 | 905 | 923 | 933 | 931 |
| 30 | 68 | 338 | 268 | 266 | 288 | 306 | 41 | 299 | 156 | 154 | 164 | 182 | 189 | 447 | 516 | 774 | 781 | 799 | 801 | 807 | 617 | 895 | 656 | 914 | 649 | 667 | 689 | 695 |
| 31 | 20 | 278 | 76 | 334 | 53 | 119 | 301 | 367 | 133 | 151 | 461 | 479 | 161 | 499 | 576 | 786 | 612 | 630 | 812 | 830 | 724 | 662 | 972 | 910 | 757 | 1015 | 685 | 943 |
| 32 | 321 | 327 | 96 | 354 | 377 | 363 | 209 | 215 | 409 | 143 | 244 | 502 | 253 | 239 | 836 | 854 | 845 | 591 | 872 | 878 | 704 | 722 | 985 | 731 | 672 | 950 | 760 | 766 |
| 33 | 12 | 10 | 36 | 102 | 196 | 194 | 168 | 234 | 520 | 586 | 544 | 550 | 644 | 710 | 744 | 750 | | | | | | | | | | | | |
| 34 | 13 | 11 | 37 | 103 | 197 | 195 | 169 | 235 | 521 | 587 | 545 | 551 | 645 | 711 | 745 | 751 | | | | | | | | | | | | |
| 35 | 8 | 14 | 100 | 38 | 192 | 198 | 232 | 170 | 584 | 522 | 548 | 546 | 708 | 646 | 748 | 746 | | | | | | | | | | | | |
| 36 | 9 | 15 | 101 | 39 | 193 | 199 | 233 | 171 | 585 | 523 | 549 | 547 | 709 | 647 | 749 | 747 | | | | | | | | | | | | |
| 37 | 284 | 282 | 308 | 374 | 468 | 466 | 440 | 506 | 792 | 858 | 816 | 822 | 916 | 982 | 1016 | 1022 | | | | | | | | | | | | |
| 38 | 285 | 283 | 309 | 375 | 469 | 467 | 441 | 507 | 793 | 859 | 817 | 823 | 917 | 983 | 1017 | 1023 | | | | | | | | | | | | |
| 39 | 280 | 286 | 372 | 310 | 464 | 470 | 504 | 442 | 856 | 794 | 820 | 818 | 980 | 918 | 1020 | 1018 | | | | | | | | | | | | |
| 40 | 281 | 287 | 373 | 311 | 465 | 471 | 505 | 443 | 857 | 795 | 821 | 819 | 981 | 919 | 1021 | 1019 | | | | | | | | | | | | |

Table 7. A covering of the 1024 vertices of $G_5$ by 40 disjoint complete subgraphs of sizes 28 and 16. The encoding is as in Table 1.

The method of getting a minimal clique covering for $G_4$ uses backtracking and the structure of the independence square (Table 3). As discussed, any clique cover for $G_4$ has at least 22 cliques. One way to search for a cover using 22 cliques is to use twenty 12-cliques and two 8-cliques. So an initial goal was to search for 20 disjoint 12-cliques. The complete set of all 86,012 12-cliques was generated. We then tried backtracking on these to find a set of 20 pairwise disjoint cliques that could extend to a clique cover but the problem size proved unmanageable. Many search paths would get stuck after including only 16 of the 12-cliques. A second problem is that if after including the twenty 12-cliques, there is a row or column in the independence square that is not covered $k$ times, then it is necessary to add at least $k$ more cliques to complete the cover. So an auspicious selection of twenty 12-cliques should leave each row uncovered at most two times each.

To attempt to deal with both of these problems, the search was restricted to only cliques that had certain subsets of the rows missing. After inspecting the subsets of rows that could be missing from one of the cliques, the following selection was made (where the rows are indexed by 0, 1,…, 15.

Group 1: Rows 0, 3, 12, and 15 are missing.
Group 2: Rows 1, 2, 13, and 14 are missing.
Group 3: Rows 4, 7, 8, and 11 are missing.
Group 4: Rows 5, 6, 9, and 10 are missing.

Backtracking on just these cliques led to several sets of twenty 12-cliques. A backtracking program was used to try to complete the clique cover, and this quickly led to a solution (most of the sets of 20 do not extend to a clique cover of size 22 but it did not take long to find one that did); see Table 6. The set of 20 that completed had 5 tuples from each group meaning that each row was used exactly 15 times (and was missing once). Similar ideas yielded the 44-clique for $G_5$ (Table 7).

It is well-known (see [DELSMW, Thm. 4.2]) that a clique in $G_d$ can be used to create a clique twice as large in $G_{d+1}$ by making two copies of it, prefacing the first copy with the digit 0, adding 1 modulo 4 to each position of each tuple in the second copy, and then prefacing each tuple in the second copy with the symbol 2. The clique also can be doubled using digits 1 and 3 as the first digits instead of 0 and 2. To double a clique cover, start with each clique $C$. Double $C$ using first digits 0 and 2. Then double the clique $C$ again using 1 and 3 as the initial digits. The original clique cover used all the $d$-tuples exactly once. For first digit 0 and 1 it is easy to see that each tuple is used exactly once. Similarly for first digits 1 and 3, each tuple appears exactly once, because adding 11…1 to every tuple of dimension $d$ gives back the complete set of tuples in dimension $d$. This construction gives a clique cover in $G_{d+1}$whose size is twice that of the cover of $G_d$. The 40-cover of $G_4$ therefore yields the upper bounds of 80 and 160 for the next two Keller graphs.

The success in finding clique covers whose size equals the lower bound suggests that this holds true in the remaining three cases, $d = 5$, 6, and 7.

**Conjecture 9.** For all $d$, $\theta(G_d) = \left\lceil \frac{4^d}{\omega(G_d)} \right\rceil$.

# 7. Conclusion

The first investigations for all our results involved computer experimentation. The patterns that one finds by such work often lead to new observations, which can sometimes be proved by classical methods. But one can be led astray. The assertion that all connected, vertex-transitive graphs (except five small examples) have Hamilton decompositions was conjectured to be true by Wagon based on extensive computations on over 100,000 graphs, including all graphs of 30 or fewer vertices. But the assertion is now known to be false [BD]. So the status of the related conjecture that all Hamiltonian vertex-transitive graphs are Hamilton-connected (except cycles and the dodecahedral graph) is a bit of a mystery. Specific conjectures such as the ones we presented here about queen graphs, Keller graphs, and the Mycielskian operation are easier to contemplate.